\documentclass[11pt]{amsart}
\usepackage{mathrsfs}

\usepackage{amsmath,amssymb,amsthm,amsfonts,graphicx}

\usepackage{url,hyperref}
\newtheorem{theorem}{Theorem}[section]
\newtheorem{proposition}[theorem]{Proposition}
\newtheorem{lemma}[theorem]{Lemma}
\newtheorem{corollary}[theorem]{Corollary}
\newtheorem{conjecture}[theorem]{Conjecture}
\theoremstyle{definition}
\newtheorem{definition}[theorem]{Definition}

\begin{document}

\title[Compact ancient solutions to three-dimensional Ricci flow]{Unique asymptotics of compact ancient solutions to three-dimensional Ricci flow}
\author[S. Angenent, S. Brendle, P. Daskalopoulos, and N. Sesum]{Sigurd Angenent, Simon Brendle, Panagiota Daskalopoulos, and Natasa Sesum}
\address{Department of Mathematics \\ University of Wisconsin \\ Madison WI 53706}
\address{Department of Mathematics \\ Columbia University \\ New York NY 10027}
\address{Department of Mathematics \\ Columbia University \\ New York NY 10027}
\address{Department of Mathematics \\ Rutgers University \\ Piscataway NJ 08854}
\begin{abstract}
We consider compact ancient solutions to the three-dimensional Ricci flow which are $\kappa$-noncollapsed. We prove that such a solutions is either a family of shrinking round spheres, or it has a unique asymptotic behavior as $t \to -\infty$ which we describe. This analysis applies in particular to the ancient solution constructed by Perelman.
\end{abstract}
\thanks{The second author was supported by the National Science Foundation under grant DMS-1806190 and by the Simons Foundation. The third author was supported by the National Science Foundation under grant DMS-1266172. The fourth author was supported by the National Science Foundation under grants DMS-1056387 and DMS-1811833.} 

\maketitle 

\tableofcontents

\section{Introduction}

A solution to a geometric evolution equation such as the Ricci flow or the Mean Curvature Flow is called \textit{ancient} if it exists for all $t\in(-\infty,T]$, for some $T$.  In the special case where the ancient solution exists for all $t \in (-\infty, \infty)$, it is called \textit{eternal}. Ancient solutions were first studied by Hamilton \cite{Hamilton3}. They play a central role in Perelman's work \cite{Perelman1},\cite{Perelman2} on singularity formation in the Ricci flow in dimension $3$. In particular, blow up limits at a singularity give rise to an ancient solution. In higher dimensions, there is a similar picture if we assume that the initial metric has positive isotropic curvature (see \cite{Hamilton4}, \cite{Brendle2}). For all these flows, the requirement that a solution should exist for all times $t \leq T$ is quite restrictive, especially when combined with a noncollapsing assumption or a positive curvature condition. In various cases, it is possible to give a complete classification; this gives a very precise description of how singularities form. 

For the two-dimensional Ricci flow, Perelman \cite{Perelman1} proved that there is only one $\kappa$-noncollapsed ancient solution, namely the shrinking spheres. Daskalopoulos, Hamilton, and \v Se\v sum \cite{Daskalopoulos-Hamilton-Sesum2} gave a complete classification of all compact ancient solutions to the Ricci flow in dimension $2$, without any noncollapsing assumptions. It turns out that the complete list contains only the shrinking spheres and the King solution. The latter was first discovered by King \cite{King} (and later independently by Rosenau \cite{Rosenau}) in the context of the logarithmic fast-diffusion equation on $\mathbb{R}^2$. The King solution also appears as the sausage model in the context of quantum field theory, in the independent work of Fateev-Onofri-Zamolodchikov \cite{Fateev-Onofri-Zamolodchikov}. While the King solution is not self-similar, it may be visualized as two cigar solitons glued together. Noncompact ancient solutions to the two-dimensional Ricci flow were classified by Daskalopoulos and \v Se\v sum \cite{Daskalopoulos-Sesum} (see also \cite{Chu}). It turns out in this case the only ancient solutions with bounded curvature are the cigar solitons. This gives a classification of all ancient solutions to the two-dimensional Ricci flow, covering both the compact and noncompact case.

Solutions analogous to the King solution exist in the higher dimensional Yamabe flow as well. Like the King solution, this is a rotationally symmetric ancient solution which is not self-similar. It can be written in closed form, and was found by King \cite{King}. However, in the case of the Yamabe flow many more ancient solutions exist. The known examples on $S^n$ include a five-parameter family of Type I ancient solutions found in \cite{Daskalopoulos-del-Pino-King-Sesum} (which includes the King solution as a special case), and the so-called "towers of bubbles" constructed in \cite{Daskalopoulos-del-Pino-Sesum} (which are of Type II). These examples suggest that it will be difficult to classify all ancient solutions to the Yamabe flow.

For curve shortening flow (i.e. mean curvature flow for curves in the plane), Daskalopoulos, Hamilton, and \v Se\v sum \cite{Daskalopoulos-Hamilton-Sesum1} classified all ancient compact convex solutions by showing that the only possibilities are the shrinking circles and the so-called Angenent ovals. In higher dimensions, White \cite{White} and Haslhofer and Hershkovits \cite{Haslhofer-Hershkovits} constructed compact ancient solutions which are rotationally symmetric but are not solitons. These can be viewed as the higher dimensional generalization of the Angenent ovals; however, no closed form expression seems to exist. For mean curvature flow in $\mathbb{R}^3$, Brendle and Choi \cite{Brendle-Choi1} classified all noncompact ancient solutions which are convex and $\alpha$-noncollapsed: the only example is the rotationally symmetric bowl soliton which moves by translations under the flow. An analogous result holds in higher dimensions, under an additional assumption that the ancient solution is uniformly two-convex (cf. \cite{Brendle-Choi2}). Angenent, Daskalopoulos, and \v Se\v sum \cite{Angenent-Daskalopoulos-Sesum2} recently classified all compact ancient solutions which are uniformly two-convex and $\alpha$-noncollapsed. They showed that, besides the shrinking spheres, there is only one example, namely the ancient oval solution constructed by White \cite{White} and by Haslhofer--Hershkovits \cite{Haslhofer-Hershkovits}. 

Our main focus in this paper will be the Ricci flow in dimension $3$. Following \cite{Perelman1}, we will consider ancient $\kappa$-solutions: these are ancient solutions which are complete; non-flat; $\kappa$-noncollapsed; and have bounded and nonnegative curvature. In the noncompact case, Perelman made the following conjecture:

\begin{conjecture}[Perelman \cite{Perelman1}]
\label{perelman.conjecture.noncompact.case}
Let $(M,g(t))$ be a noncompact ancient $\kappa$-solution to the Ricci flow in dimension $3$ with positive curvature. Then $(M,g(t))$ is isometric to the Bryant soliton, up to scaling.
\end{conjecture}

The first major step towards Perelman's conjecture was carried out in \cite{Brendle1}, where it was shown that the Bryant soliton is the only steady gradient Ricci soliton in dimension $3$ which is $\kappa$-noncollapsed and has positive curvature. Perelman's conjecture was recently solved in full generality in \cite{Brendle3}. The proof in \cite{Brendle3} consists of two parts. In the first part, it is shown that any noncompact ancient $\kappa$-solution which is rotationally symmetric must be the Bryant soliton. In the second part, it is shown that every noncompact ancient $\kappa$-solution is, in fact, rotationally symmetric. 

We now turn to the compact case. A compact ancient $\kappa$-solution in dimension $3$ must have positive curvature; in particular, its universal cover must be diffeomorphic to $S^3$. Perelman established the existence of a rotationally symmetric ancient solution on $S^3$. This ancient solution is of Type II, i.e. $\limsup_{t \to -\infty} (-t) \, R_{\text{\rm max}}(t) = \infty$. Perelman's solution can be viewed as the three-dimensional analogue of the two-dimensional King solution. However, unlike the King solution (which is collapsed), Perelman's ancient solution is $\kappa$-noncollapsed. Going forward in time, Perelman's ancient solution shrinks to a round point. As $t \to -\infty$, Perelman's ancient solution looks like two Bryant solitons glued together.

The following conjecture can be viewed as the analogue of Perelman's Conjecture \ref{perelman.conjecture.noncompact.case} in the compact setting:

\begin{conjecture}
\label{perelman.conjecture.compact.case}
Let $(S^3,g(t))$ be a compact ancient $\kappa$-solution to the Ricci flow in dimension $3$. Then $g(t)$ is either a family of contracting spheres or Perelman's ancient solution.
\end{conjecture}

As pointed out in \cite{Brendle3}, the methods in that paper imply that any compact ancient $\kappa$-solution in dimension $3$ must be rotationally symmetric. The classification of compact ancient solutions with rotational symmetry is a difficult problem. A major challenge in this problem comes from the fact that Perelman's solution is not given in explicit form and is not a soliton. A similar challenge appears in the classification of compact ancient solutions to mean curvature flow which was resolved in \cite{Angenent-Daskalopoulos-Sesum2}. To overcome this problem, one needs a very precise understanding of the asymptotic behavior of the ancient solution as $t \to -\infty$. In this paper, we carry out the necessary asymptotic analysis for compact rotationally symmetric ancient solutions to the three-dimensional Ricci flow:

\begin{theorem}
\label{asymptotics.compact.case}
Let $(S^3,g(t))$ be a rotationally symmetric ancient $\kappa$-solution which is not isometric to a family of shrinking spheres. Then we can find a reference point $q \in S^3$ such that the following holds. Let $F(z,t)$ denote the radius of the sphere of symmetry in $(S^3,g(t))$ which has signed distance $z$ from the reference point $q$. Then the profile $F(z,t)$ has the following asymptotic expansions:
\begin{itemize}
\item[(i)] Fix a large number $L$. Then, as $t \to -\infty$, we have 
\[F(z,t)^2 = -2t - \frac{z^2+2t}{2\log(-t)} + o \Big ( \frac{(-t)}{\log(-t)} \Big )\] 
for $|z| \leq L\sqrt{-t}$
\item[(ii)] Fix a small number $\theta>0$. Then as $t \to -\infty$, we have 
\[F(z,t)^2 = -2t - \frac{z^2}{2\log(-t)} + o(-t)\] 
for $|z| \leq 2\sqrt{1-\theta^2} \sqrt{(-t) \log(-t)}$.
\item[(iii)] The reference point $q$ has distance $(2+o(1)) \sqrt{(-t) \log(-t)}$ from each tip. The scalar curvature at each tip is given by $(1+o(1)) \, \frac{\log(-t)}{(-t)}$. Finally, if we rescale the solution around one of the tips, then the rescaled solutions converge to the Bryant soliton as $t \to -\infty$.
\end{itemize}
\end{theorem}

In a recent work \cite{Brendle-Daskalopoulos-Sesum}, we use Theorem \ref{asymptotics.compact.case} to settle Conjecture \ref{perelman.conjecture.compact.case}, in a similar way that results about unique asymptotics of ancient ovals shown in \cite{Angenent-Daskalopoulos-Sesum1} were used to prove the classification result of closed ancient mean curvature flow solutions (see \cite{Angenent-Daskalopoulos-Sesum2}). The proof in \cite{Brendle-Daskalopoulos-Sesum} turned out to be very involved and hence was written in a separate paper.

In order to prove Theorem \ref{asymptotics.compact.case} we will combine techniques developed in \cite{Angenent-Daskalopoulos-Sesum1} and \cite{Brendle3}. In \cite{Brendle3}, under the assumption on rotational symmetry, Brendle constructed barriers by using gradient Ricci solitons with singularity at the tip which were found by Bryant \cite{Bryant}. In the proof of Theorem \ref{asymptotics.compact.case}, we use these barriers to localize our equation in the cylindrical region (where the solution is close to a cylinder of radius $\sqrt{-2t}$). Similar localization arguments were employed in \cite{Angenent-Daskalopoulos-Sesum1} and \cite{Brendle3}. The localization enables us to do spectral decomposition in the cylindrical region and obtain refined asymptotics of our solution in the cylindrical region.

The outline of the paper is as follows: In Section \ref{basic.properties}, we record some basic properties of compact ancient $\kappa$-solutions. In particular, we show that, if $-t$ is sufficiently large, the solution looks like the Bryant soliton near each tip. In Section \ref{asymptotics.near.cylinder}, we use the barriers from \cite{Brendle3} to achieve the spectral decomposition of our solution. This allows us to apply the Merle-Zaag lemma (see \cite{Merle-Zaag}). This leaves us with two possibilities: either the positive modes dominate, or the neutral mode dominates. The former case is ruled out in Section \ref{positive.modes.cannot.dominate}. In the latter case, we obtain precise asymptotics in the cylindrical region (see Section \ref{neutral.mode.dominates}). Subsequently, we combine this exact behavior in the cylindrical region with barrier arguments to obtain the precise behavior of our solution in the intermediate region (see Section \ref{intermediate.region.asymptotics}).  In Section \ref{tip.region.asymptotics}, we obtain the precise behavior of the distance from the reference point $q$ to each tip. Combining this estimate with Hamilton's Harnack inequality \cite{Hamilton2}, we obtain precise asymptotics for the scalar curvature at each tip. Finally, in Appendix \ref{estimate.for.1d.heat.equation}, we state and prove an elementary estimate for the one-dimensional heat equation, which is needed in Section \ref{positive.modes.cannot.dominate} in the proof of Proposition \ref{key.estimate}. In Appendix \ref{properties.of.Bryant.soliton} we state known facts about the Bryant Ricci soliton. \\

\section{Basic properties of compact ancient solutions}

\label{basic.properties}

Throughout this paper, we assume that $(S^3,g(t))$, $t \in (-\infty,0]$, is an ancient $\kappa$-solution which is rotationally symmetric. Moreover, we assume that $(S^3,g(t))$ is not a family of shrinking round spheres.

\begin{lemma}
\label{asymptotic.soliton}
The asymptotic soliton of $(S^3,g(t))$ is a cylinder. In other words, suppose that we fix a point $q \in S^3$. Consider a sequence of times $t_k \to -\infty$ and a sequence of points $x_k \in S^3$ such that $\sup_k \ell(x_k,t_k) < \infty$, where $\ell$ denotes the reduced distance from $(q,0)$. If we dilate the manifold $(S^3,g(t_k))$ around the point $x_k$ by the factor $(-t_k)^{-\frac{1}{2}}$, then the rescaled manifolds converge to a cylinder of radius $\sqrt{2}$.
\end{lemma}

\textbf{Proof.} 
By work of Perelman \cite{Perelman1}, the rescaled manifolds converge in the Cheeger-Gromov sense to a smooth limit, and the limit is a shrinking gradient Ricci soliton. We claim that the limiting soliton must be a cylinder of radius $\sqrt{2}$. We distinguish two cases:

\textit{Case 1:} If the limiting soliton is compact, then it must have constant sectional curvature $\frac{1}{4}$. In particular, the sectional curvatures of $(S^3,g(t_k))$ lie in the interval $[\frac{1-\varepsilon_k}{(-4t_k)},\frac{1+\varepsilon_k}{(-4t_k)}]$, where $\varepsilon_k \to 0$ as $k \to \infty$. Hamilton's curvature pinching estimates now imply that $(S^3,g(t))$ has constant sectional curvature for each $t$ (cf. \cite{Hamilton1}). Thus, $(S^3,g(t))$ is a family of shrinking round spheres, contrary to our assumption.

\textit{Case 2:} If the limiting soliton is noncompact, then results of Perelman imply that the limit is a cylinder of radius $\sqrt{2}$ (cf. \cite{Perelman2}, Section 1.1). This proves the assertion. \\

\begin{lemma}
\label{diameter}
Given any sequence of times $t_k \to -\infty$, we have 
\[R_{\text{\rm max}}(t_k) \, \text{\rm diam}_{g(t_k)}(S^3,g(t_k))^2 \to \infty,\] 
where $R$ denotes the scalar curvature and $R_{\text{\rm max}}$ denotes the maximum of the scalar curvature over all points in space.
\end{lemma}

\textbf{Proof.} 
By work of Perelman \cite{Perelman1}, we can find a sequence of points $x_k \in S^3$ such that $\ell(x_k,t_k) \leq 3$ for each $k$. By Lemma \ref{asymptotic.soliton}, if we dilate the manifold $(S^3,g(t_k))$ around the point $x_k$ by the factor $(-t_k)^{-\frac{1}{2}}$, then the rescaled manifolds converge to a cylinder of radius $\sqrt{2}$. From this, the assertion follows. \\

\begin{proposition}
\label{possible.limits}
Consider a sequence of times $t_k \to -\infty$ and an arbitrary sequence of points $x_k \in S^3$. If we dilate the flow around the point $(x_k,t_k)$ by the factor $R(x_k,t_k)^{\frac{1}{2}}$, then (after passing to a subsequence) the rescaled flows converge to either a family of shrinking cylinders or the Bryant soliton.
\end{proposition}

\textbf{Proof.} By Perelman's compactness theorem for ancient $\kappa$-solutions, the rescaled flows converge to an ancient $\kappa$-solution. If the limit is compact, then $\sup_k R_{\text{\rm max}}(t_k) \, \text{\rm diam}_{g(t_k)}(S^3,g(t_k)) < \infty$, which contradicts Lemma \ref{diameter}. Consequently, the limit must be noncompact. The results in \cite{Brendle3} now imply that the limit is either a family of shrinking cylinders or the Bryant soliton. \\

\begin{proposition}
\label{rescaling.around.tip}
Consider a sequence of times $t_k \to -\infty$. If we rescale the solution around one of the tips, then the rescaled solutions converge to the Bryant solution.
\end{proposition}

\textbf{Proof.} 
By symmetry, the tracefree part of the Ricci tensor vanishes at the tip. Consequently, if we rescale around the tip, the limit cannot be a cylinder. By Proposition \ref{possible.limits}, the only possible limit is the Bryant soliton. \\

\begin{corollary}
\label{scalar.curvature.at.tips.1}
Let $R_{\text{\rm tip},1}(t)$ and  $R_{\text{\rm tip},2}(t)$ denote the scalar curvature at the tips. Then $\frac{d}{dt} R_{\text{\rm tip},1}(t) \leq o(1) \, R_{\text{\rm tip},1}(t)^2$ and $\frac{d}{dt} R_{\text{\rm tip},2}(t) \leq o(1) \, R_{\text{\rm tip},2}(t)^2$.
\end{corollary}

\textbf{Proof.} This is a direct consequence of Proposition \ref{rescaling.around.tip}. \\

\begin{corollary}
\label{scalar.curvature.at.tips.2} 
We have $(-t) \, R_{\text{\rm tip},1}(t) \to \infty$ and $(-t) \, R_{\text{\rm tip},2}(t) \to \infty$ as $t \to -\infty$. 
\end{corollary}

\textbf{Proof.} 
This follows by integrating the differential inequality in Corollary \ref{scalar.curvature.at.tips.1}. \\

\begin{definition}
\label{evolving.neck}
Let $(x_0,t_0)$ be a point in space-time with $R(x_0,t_0) = r^{-2}$. We say that $(x_0,t_0)$ lies at the center of an evolving $\varepsilon$-neck if, after rescaling by the factor $r^{-1}$, the parabolic neighborhood $B_{g(t_0)}(x_0,\varepsilon^{-1} r) \times [t_0-\varepsilon^{-1} r^2,t_0]$ is $\varepsilon$-close in $C^{[\varepsilon^{-1}]}$ to a family of shrinking cylinders.
\end{definition}

\begin{proposition}
\label{neck}
Given $\varepsilon>0$ and $\delta>0$, we can find a time $t_0 = t_0(\varepsilon,\delta)$ so that the following holds. Suppose $(x,t)$ is a point in spacetime such that $t \leq t_0$, and the radius of the sphere of symmetry through $(x,t)$ is at least $\delta \sqrt{-2t}$. Then $(x,t)$ lies at the center of an evolving $\varepsilon$-neck.
\end{proposition}

\textbf{Proof.} 
We argue by contradiction. Suppose that there exists a sequence of points $(x_k,t_k)$ in space-time with the following properties: 
\begin{itemize}
\item $t_k \to -\infty$.
\item The sphere of symmetry through $(x_k,t_k)$ has radius at least $\delta \sqrt{-2t_k}$. 
\item The point $(x_k,t_k)$ does not lie at the center of an evolving $\varepsilon$-neck.
\end{itemize}
By assumption, the sphere of symmetry through $(x_k,t_k)$ has radius $r_k \geq \delta \sqrt{-2t_k}$. At the point $(x_k,t_k)$, the sectional curvature of the plane tangent to the sphere of symmetry is bounded from above by $r_k^{-2}$. Consequently, the minimum sectional curvature at $(x_k,t_k)$ satisfies $K_{\text{\rm min}}(x_k,t_k) \leq r_k^{-2} \leq \frac{1}{\delta^2 (-2t_k)}$.

Since the point $(x_k,t_k)$ does not lie at the center of an evolving $\varepsilon$-neck, we must have $\liminf_{k \to \infty} R(x_k,t_k)^{-1} \, K_{\text{\rm min}}(x_k,t_k) > 0$. Since $(-t_k) \, K_{\text{\rm min}}(x_k,t_k) \leq \frac{1}{2\delta^2}$, it follows that $\limsup_{k \to \infty} (-t_k) \, R(x_k,t_k) < \infty$. On the other hand, we have $(-t_k) \, R_{\text{\rm tip},1}(t_k) \to \infty$ and $(-t_k) \, R_{\text{\rm tip},2}(t_k) \to \infty$ by Corollary \ref{scalar.curvature.at.tips.2}. This gives $R(x_k,t_k)^{-1} \, R_{\text{\rm tip},1}(t_k) \to \infty$ and $R(x_k,t_k)^{-1} \, R_{\text{\rm tip},2}(t_k) \to \infty$. By Perelman's long range curvature estimate, the distance of $(x_k,t_k)$ from each tip is bounded from below $A_k \, R(x_k,t_k)^{-\frac{1}{2}}$, where $A_k \to \infty$. Hence, if we dilate the flow around the point $(x_k,t_k)$ by the factor $R(x_k,t_k)^{\frac{1}{2}}$ and pass to the limit, then the limit contains a line. By the Cheeger-Gromov splitting theorem, the limit splits as a product. Thus, the limit is a cylinder. This contradicts the fact that $\liminf_{k \to \infty} R(x_k,t_k)^{-1} \, K_{\text{\rm min}}(x_k,t_k) > 0$. This is a contradiction. \\

\section{Asymptotic analysis near the cylinder}

\label{asymptotics.near.cylinder}

We continue to assume that $(S^3,g(t))$, $t \in (-\infty,0]$, is an ancient $\kappa$-solution which is rotationally symmetric. Moreover, we assume that $(S^3,g(t))$ is not a family of shrinking round spheres. We begin by fixing a base point $q$. This point will be chosen such that $\limsup_{t \to -\infty} (-t) \, R(q,t) \leq 100$. The existence of such a point follows from the Neck Stability Theorem of Kleiner and Lott \cite{Kleiner-Lott}. The result in \cite{Kleiner-Lott} is stated in the noncompact setting, but the argument can be easily adapted to the compact case:

\begin{proposition}[cf. Kleiner-Lott \cite{Kleiner-Lott}, Section 6]
\label{neck.stability}
There exists a point $q \in S^3$ with the property that $\limsup_{t \to -\infty} (-t) \, R(q,t) \leq 100$. 
\end{proposition}

\textbf{Proof.} 
Suppose that the assertion is false. By Lemma \ref{asymptotic.soliton}, we can find a sequence of points $(q_k,s_k)$ in space-time and a sequence of positive numbers $\varepsilon_k \to 0$ with the property that $(q_k,s_k)$ lies at the center of an evolving $\varepsilon_k$-neck. Our assumption implies $\limsup_{t \to -\infty} (-t) \, R(q_k,t) > 100$ for each $k$. This implies $\limsup_{t \to -\infty} (s_k-t) \, R(q_k,t) > 100$ for each $k$. For each $k$, we define $t_k := \sup \{t \leq s_k: (s_k-t) \, R(q_k,t) > 10\}$. Clearly, $t_k \leq s_k - \varepsilon_k^{-1} \, R(q_k,s_k)^{-1}$, $(s_k-t_k) \, R(q_k,t_k) = 10$, and $(s_k-t) \, R(q_k,t) \leq 10$ for all $t \in [t_k,s_k]$. 

Let $\ell_k(x,t)$ denote the reduced distance of $(x,t)$ from $(q_k,s_k)$, and let $V_k(t) = \int (s_k-t)^{-\frac{3}{2}} \, e^{-\ell_k(x,t)} \, d\text{\rm vol}_{g(t)}$ denote the reduced volume. By definition of $t_k$, we have $\ell(q_k,t_k) \leq \frac{1}{2\sqrt{s_k-t_k}} \int_{t_k}^{s_k} \sqrt{s_k-t} \, R(q_k,t) \, dt \leq 100$. Since $(q_k,s_k)$ lies on an evolving $\varepsilon_k$-neck, we obtain $\limsup_{k \to \infty} V_k(s_k + \tau \, R(q_k,s_k)^{-1}) \leq V_{\text{\rm cyl}}(\tau)$ for each $\tau \in (-\infty,0)$, where $V_{\text{\rm cyl}}(\tau)$ denotes the reduced volume of a family of shrinking cylinders. This implies $\limsup_{k \to \infty} V_k(s_k - \varepsilon_k^{-1} \, R(q_k,s_k)^{-1}) \leq V_{\text{\rm cyl}}(-\infty)$. On the other hand, since the asymptotic soliton of $(S^3,g(t))$ is a cylinder, Perelman's monotonicity formula gives $V_k(t) \geq V_{\text{\rm cyl}}(-\infty)$ for all $k$ and all $t \in (-\infty,s_k)$.

We now dilate the flow around the point $(q_k,t_k)$ by the factor $(s_k-t_k)^{-\frac{1}{2}}$. By work of Perelman, the rescaled flows converge in the Cheeger-Gromov sense to a smooth limit. On the limit, the reduced volume is constant, and equals $V_{\text{\rm cyl}}(-\infty)$. Consequently, the limit must be a shrinking gradient Ricci  soliton. Since $(s_k-t_k) \, R(q_k,t_k) = 10$ for each $k$, the limit is non-flat. Moreover, the limit cannot have constant curvature, for otherwise our ancient solution $(S^3,g(t))$ would have constant curvature, contradicting our assumption. Thus, the limiting gradient soliton must be a cylinder with scalar curvature $1$. In particular, $(s_k-t_k) \, R(q_k,t_k) \to 1$ as $k \to \infty$. This contradicts the fact that $(s_k-t_k) \, R(q_k,t_k) = 10$ for each $k$. \\

\begin{proposition}
\label{limit.around.q}
Consider a sequence of times $t_k \to -\infty$. If we dilate the flow around the point $(q,t_k)$ by the factor $(-t_k)^{-\frac{1}{2}}$, then the rescaled manifolds converge to a cylinder of radius $\sqrt{2}$.
\end{proposition}

\textbf{Proof.} 
By our choice of $q$, we have $\limsup_{t \to -\infty} (-t) \, R(q,t) < \infty$. Consequently, $\lim_{t \to -\infty} \ell(q,t) < \infty$, where $\ell$ denotes the reduced distance from $(q,0)$. Hence, the assertion follows from Lemma \ref{asymptotic.soliton}. \\

For each $t$, we denote by $F(z,t)$ the radius of a sphere of symmetry which has signed distance $z$ from the reference point $q$. In particular, $F(0,t)$ is the radius of the sphere of symmetry passing through $q$. Since the manifold has positive sectional curvature, we have $F_{zz} \leq 0$. Moreover, we have $F_z = \pm 1$ at the tips. Consequently, $-1 \leq F_z \leq 1$ at each point in space-time. 

\begin{proposition}
\label{evolution.of.F}
The function $F$ satisfies the evolution equation 
\begin{align*} 
F_t(z,t)  
&= F_{zz}(z,t) - F(z,t)^{-1} \, (1-F_z(z,t)^2) \\ 
&- 2 \, F_z(z,t) \int_0^z \frac{F_{zz}(z',t)}{F(z',t)} \, dz' \\ 
&= F_{zz}(z,t) - F(z,t)^{-1} \, (1+F_z(z,t)^2) \\ 
&+ 2 \, F_z(z,t) \, \bigg [ F(0,t)^{-1} \, F_z(0,t) - \int_0^z \frac{F_z(z',t)^2}{F(z',t)^2} \, dz' \bigg ]. 
\end{align*} 
\end{proposition}

\textbf{Proof.} 
Let us consider an arbitrary point $p \in S^3$. Let $r(t)$ denote the radius (with respect to the metric $g(t)$) of the sphere of symmetry passing through the point $p$, and let $z(t)$ denote the signed distance of that sphere from the reference point $q$. Clearly, $r(t) = F(z(t),t)$, hence 
\[r'(t) = F_t(z(t),t) + F_z(z(t),t) \, z'(t).\] 
Using the formula 
\[\text{\rm Ric} = -\frac{F_{zz}}{F} \, (g + dz \otimes dz) + \frac{1-F_z^2}{F^2} \, (g - dz \otimes dz),\] 
we obtain 
\[z'(t) = 2 \int_0^{z(t)} \frac{F_{zz}(z',t)}{F(z',t)} \, dz'\] 
and 
\[r'(t) = F_{zz}(z(t),t) - F(z(t),t)^{-1} \, (1-F_z(z(t),t)^2).\] 
Putting these facts together, we obtain 
\begin{align*} 
F_t(z,t) 
&= F_{zz}(z,t) - F(z,t)^{-1} \, (1-F_z(z,t)^2) \\ 
&- 2 \, F_z(z,t) \int_0^z \frac{F_{zz}(z',t)}{F(z',t)} \, dz'. 
\end{align*}
Integration by parts gives 
\begin{align*} 
F_t(z,t) 
&= F_{zz}(z,t) - F(z,t)^{-1} \, (1+F_z(z,t)^2) \\ 
&+ 2 \, F_z(z,t) \, \bigg [ F(0,t)^{-1} \, F_z(0,t) - \int_0^z \frac{F_z(z',t)^2}{F(z',t)^2} \, dz' \bigg ]. 
\end{align*}
This completes the proof. \\

\begin{definition}
\label{def.r_max}
For each $t$, let $r_{\text{\rm max}}(t) = \sup_z F(z,t)$ denote the maximum radius at time $t$. 
\end{definition}

\begin{lemma}
\label{lower.bound.for.r_max}
We have $-\frac{1}{2} \, \frac{d}{dt} (r_{\text{\rm max}}(t)^2) \geq 1$. In particular, $r_{\text{\rm max}}(t) \geq \sqrt{-2t}$ for each $t$. 
\end{lemma}

\textbf{Proof.} 
Consider the point where the radius is maximal. At that point, $F = r_{\text{\rm max}}(t)$, $F_z = 0$, and $F_{zz} \leq 0$. Using the evolution equation for $F$, we conclude that $-\frac{1}{2} \, \frac{\partial}{\partial t} (F^2) \geq 1$ at the point where the radius is maximal. Thus, we conclude that $-\frac{1}{2} \, \frac{d}{dt} (r_{\text{\rm max}}(t)^2) \geq 1$. Integrating over $t$, we obtain $r_{\text{\rm max}}(t) \geq \sqrt{-2t}$ for each $t$. This completes the proof of Lemma \ref{lower.bound.for.r_max}. \\

\begin{lemma}
\label{upper.bound.for.r_max}
We have $-\frac{1}{2} \, \frac{d}{dt} (r_{\text{\rm max}}(t)^2) \leq 1+o(1)$ as $t \to -\infty$. In particular, $r_{\text{\rm max}}(t) \leq (1+o(1)) \, \sqrt{-2t}$ as $t \to -\infty$. 
\end{lemma}

\textbf{Proof.} 
Let $\varepsilon > 0$ be given. Let us consider the point where the radius is maximal. At that point, $F = r_{\text{\rm max}}(t) \geq \sqrt{-2t}$ and $F_z = 0$. Proposition \ref{neck} implies that the point where the radius is maximal lies on an evolving $\varepsilon$-neck if $-t$ is sufficiently large. Hence, if $-t$ is sufficiently large, then we have $F_{zz} \geq -\varepsilon \, F^{-1}$ at the point where the radius is maximal. Using the evolution equation for $F$, we obtain $-\frac{1}{2} \, \frac{\partial}{\partial t} (F^2) \leq 1+\varepsilon$ at the point where the radius is maximal. Thus, we conclude that $-\frac{1}{2} \, \frac{d}{dt} (r_{\text{\rm max}}(t)^2) \leq 1+\varepsilon$ if $-t$ is sufficiently large. Since $\varepsilon > 0$ is arbitrary, it follows that $-\frac{1}{2} \, \frac{d}{dt} (r_{\text{\rm max}}(t)^2) \leq 1+o(1)$ as $t \to -\infty$. This finally implies $r_{\text{\rm max}}(t) \leq (1+o(1)) \, \sqrt{-2t}$ as $t \to -\infty$. This completes the proof of Lemma \ref{upper.bound.for.r_max}. \\

We next describe a barrier argument which plays a key role in our analysis. We begin by recalling a result from \cite{Brendle3}:

\begin{proposition}[cf. \cite{Brendle3}, Section 2]
\label{properties.of.psi_a}
There exists a positive real number $r_*$ with the following property. Given any large number $a$, we can find a continuously differentiable function $\psi_a: [r_* a^{-1},1+\frac{1}{100} \, a^{-2}] \to \mathbb{R}$ such that 
\[\psi_a(s) \psi_a''(s) - \frac{1}{2} \, \psi_a'(s)^2 + s^{-2} \, (1-\psi_a(s)) \, (s \psi_a'(s) + 2 \psi_a(s)) - s \psi_a'(s) < 0\] 
for all $s \in [r_* a^{-1},1+\frac{1}{100} \, a^{-2}]$. The function $\psi_a$ satisfies $\psi_a(s) \leq Ca^{-2}$ for all $s \in [\frac{1}{10},1+\frac{1}{100} \, a^{-2}]$. Moreover, $\psi_a(s) \geq \frac{1}{32} \, a^{-4}$ for all $s \in [r_* a^{-1},1+\frac{1}{100} \, a^{-2}]$, and $\psi_a(s) \geq a^{-2} (s^{-2}-1) + \frac{1}{16} \, a^{-4}$ for all $s \in [1-\theta,1+\frac{1}{100} \, a^{-2}]$, where $\theta$ is a small positive number. Finally, $\psi_a(r_* a^{-1}) \geq \frac{3}{2}$.
\end{proposition}

\textbf{Proof.} 
It follows from Definition 2.6 in \cite{Brendle3} that $\psi_a(s) \leq Ca^{-2}$ for all $s \in [\frac{1}{10},1+\frac{1}{100} \, a^{-2}]$ and $\psi_a(r_* a^{-1}) = 2+O(a^{-1})$. In particular, $\psi_a(r_* a^{-1}) \geq \frac{3}{2}$ if $a$ is sufficiently large. All the remaining statements follow from Proposition 2.7 in \cite{Brendle3}. \\

\begin{proposition}
\label{barrier}
Suppose that $a$ is sufficiently large. Moreover, suppose that $\bar{t}$ is a time with the property that $\frac{r_{\text{\rm max}}(t)}{\sqrt{-2t}} \leq 1+\frac{1}{100} \, a^{-2}$ for all $t \leq \bar{t}$. Then $F_z(z,t)^2 < \psi_a \big ( \frac{F(z,t)}{\sqrt{-2t}} \big )$ whenever $t \leq \bar{t}$ and $F(z,t) \geq r_* a^{-1} \sqrt{-2t}$.
\end{proposition}

\textbf{Proof.} 
Our assumption implies that $\frac{F(z,t)}{\sqrt{-2t}} \leq \frac{r_{\text{\rm max}}(t)}{\sqrt{-2t}} \leq 1+\frac{1}{100} \, a^{-2}$ for all $t \leq \bar{t}$ and all $z$. Let $\mathcal{I}$ denote the set of all times $t \leq \bar{t}$ with the property that $F_z(z,t)^2 < \psi_a \big ( \frac{F(z,t)}{\sqrt{-2t}} \big )$ whenever $F(z,t) \geq r_* a^{-1} \sqrt{-2t}$. It follows from Proposition \ref{neck} that $t \in \mathcal{I}$ if $-t$ is sufficiently large (depending on $a$). 

We claim that $\mathcal{I} = (-\infty,\bar{t}]$. Suppose this is false. Let $t_0$ denote the infimum of the set $(-\infty,\bar{t}] \setminus \mathcal{I}$. We can find a point $z_0$ such that $F(z_0,t_0) \geq r_* a^{-1} \sqrt{-2t_0}$ and $F_z(z_0,t_0)^2 = \psi_a \big ( \frac{F(z_0,t_0)}{\sqrt{-2t_0}} \big )$. Clearly, $F_z(z_0,t_0) \neq 0$ since $\psi_a$ is positive; in other words, the function $z \mapsto F(z,t_0)$ does not attain its maximum at $z_0$.

If $F(z_0,t_0) = r_* a^{-1} \sqrt{-2t_0}$, then 
\[F_z(z_0,t_0)^2 \leq 1 < \psi_a(r_* a^{-1}) = \psi_a \Big ( \frac{F(z_0,t_0)}{\sqrt{-2t_0}} \Big ).\] 
This contradicts our choice of $(z_0,t_0)$. Thus, we conclude that $F(z_0,t_0) > r_* a^{-1} \sqrt{-2t_0}$.

Since $F_z(z_0,t_0) \neq 0$, we can find a smooth function $u(r,t)$ such that $F_z(z,t)^2 = u(F(z,t),t)$ in a neighborhood of the point $(z_0,t_0)$. The function $u$ is defined locally in a neighborhood of $(r_0,t_0)$, where $r_0 := F(z_0,t_0)$. Clearly, $u(r_0,t_0) = \psi_a \big ( \frac{r_0}{\sqrt{-2t_0}} \big )$. Moreover, since $(-\infty,t_0) \subset \mathcal{I}$, we obtain $u(r,t) < \psi_a \big ( \frac{r}{\sqrt{-2t}} \big )$ for $t < t_0$. As in \cite{Brendle3}, the function $u(r,t)$ satisfies 
\[u_t = uu_{rr} - \frac{1}{2} \, u_r^2 + r^{-2} \, (1-u) \, (ru_r+2u).\] 
On the other hand, the function $\Psi_a(r,t) := \psi_a \big ( \frac{r}{\sqrt{-2t}} \big )$ satisfies 
\[\Psi_{a,t} > \Psi_a \Psi_{a,rr} - \frac{1}{2} \, \Psi_{a,r}^2 + r^{-2} \, (1-\Psi_a) \, (r \Psi_{a,r}+2\Psi_a).\] 
This contradicts the parabolic maximum principle. This completes the proof of Proposition \ref{barrier}. \\

We now perform a rescaling. We define 
\[G(\xi,\tau) := e^{\frac{\tau}{2}} \, F(e^{-\frac{\tau}{2}} \xi,-e^{-\tau}) - \sqrt{2}.\] 
Then 
\begin{align*} 
G_\tau(\xi,\tau) 
&= G_{\xi\xi}(\xi,\tau) - \frac{1}{2} \, \xi \, G_\xi(\xi,\tau) \\
&+ \frac{1}{2} \, (\sqrt{2}+G(\xi,\tau)) - (\sqrt{2}+G(\xi,\tau))^{-1} \, (1+G_\xi(\xi,\tau)^2) \\ 
&+ 2 \, G_\xi(\xi,\tau) \, \bigg [ (\sqrt{2}+G(0,\tau))^{-1} \, G_\xi(0,\tau) - \int_0^\xi \frac{G_\xi(\xi',\tau)^2}{(\sqrt{2}+G(\xi',\tau))^2} \, d\xi' \bigg ]. 
\end{align*} 
Proposition \ref{limit.around.q} implies that, as $\tau \to -\infty$, the functions $G(\xi,\tau)$ converges to $0$ in $C_{\text{\rm loc}}^\infty$. \\

\begin{definition}
\label{def.rho_max}
For each $\tau$, let $\rho_{\text{\rm max}}(\tau) := \sup_\xi G(\xi,\tau) = e^{\frac{\tau}{2}} \, r_{\text{\rm max}}(-e^{-\tau}) - \sqrt{2}$. 
\end{definition}

\begin{lemma}
\label{rho_max}
We have $\rho_{\text{\rm max}}(\tau) \geq 0$ for each $\tau$. Moreover, $\rho_{\text{\rm max}}(\tau) \to 0$ as $\tau \to -\infty$. 
\end{lemma}

\textbf{Proof.}
This is an immediate consequence of Lemma \ref{lower.bound.for.r_max} and Lemma \ref{upper.bound.for.r_max}. \\

For each time $\bar{\tau}$, we define 
\[\delta(\bar{\tau}) := \sup_{\tau \leq \bar{\tau}} (|G(0,\tau)| +\rho_{\text{\rm max}}(\tau)).\]
By definition, $\delta(\bar{\tau})$ is an increasing function of $\bar{\tau}$. Moreover, $\delta(\bar{\tau}) \to 0$ as $\bar{\tau} \to -\infty$.

\begin{lemma}
\label{derivative.of.delta}
For $-\bar{\tau}$ sufficiently large, the function $\bar{\tau} \mapsto \delta(\bar{\tau})$ is Lipschitz continuous with Lipschitz constant $1$. In particular, the function $\bar{\tau} \mapsto \delta(\bar{\tau})$ is differentiable almost everywhere, and $0 \leq \delta'(\bar{\tau}) \leq 1$ for $-\bar{\tau}$ sufficiently large. 
\end{lemma}

\textbf{Proof.} 
If $-\tau$ is sufficiently small, then the functions $\tau \mapsto G(0,\tau)$ and $\tau \mapsto G_\xi(0,\tau)$ are Lipschitz continuous with Lipschitz constant $\frac{1}{4}$. Moreover, we have shown above that $\frac{1}{2} \, \frac{d}{dt} (r_{\text{\rm max}}(t)^2) \to -1$ as $t \to -\infty$. This implies $\frac{d}{d\tau} \rho_{\text{\rm max}}(\tau) \to 0$ as $\tau \to -\infty$. Hence, if $-\tau$ is sufficiently large, then the function $\tau \mapsto \rho_{\text{\rm max}}(\tau)$ is Lipschitz continuous with Lipschitz constant $\frac{1}{4}$. Putting these facts together, we conclude that the function $\bar{\tau} \mapsto \delta(\bar{\tau})$ is Lipschitz continuous with Lipschitz constant $1$. This proves the first statement. The second statement follows from Rademacher's theorem. \\

\begin{proposition}
\label{barrier.rescaled.picture}
Fix $\bar{\tau}$, and let $a := \frac{1}{10} \, \delta(\bar{\tau})^{-\frac{1}{2}}$. Then we have $G_\xi(\xi,\tau)^2 \leq \psi_a \big ( 1+\frac{G(\xi,\tau)}{\sqrt{2}} \big )$ whenever $\tau \leq \bar{\tau}$ and $G(\xi,\tau) \geq (r_* a^{-1} - 1) \sqrt{2}$, where $r_*$ is defined as in \cite{Brendle3}.
\end{proposition}

\textbf{Proof.} 
By definition of $\delta(\bar{\tau})$, we have $\rho_{\text{\rm max}}(\tau) \leq \delta(\bar{\tau})$ for all $\tau \leq \bar{\tau}$. This implies $\frac{r_{\text{\rm max}}(t)}{\sqrt{-2t}} \leq 1+\delta(\bar{\tau})$ for all $t \leq -e^{-\bar{\tau}}$. Hence, we may apply Proposition \ref{barrier} with $\bar{t} = -e^{-\bar{\tau}}$ and $a = \frac{1}{10} \, \delta(\bar{\tau})^{-\frac{1}{2}}$. Using Proposition \ref{barrier}, we conclude that $F_z(z,t)^2 < \psi_a \big ( \frac{F(z,t)}{\sqrt{-2t}} \big )$ whenever $t \leq \bar{t}$ and $F(z,t) \geq r_* a^{-1} \sqrt{-2t}$. In other words, $G_\xi(\xi,\tau)^2 < \psi_a \big ( 1+\frac{G(\xi,\tau)}{\sqrt{2}} \big )$ whenever $\tau \leq \bar{\tau}$ and $G(\xi,\tau) \geq (r_* a^{-1} - 1) \sqrt{2}$. \\

\begin{lemma}
\label{C1.bound.for.G}
We have $|G(\xi,\tau)|+|G_\xi(\xi,\tau)| \leq C \, \delta(\tau)^{\frac{1}{4}}$ for $|\xi| \leq 2\delta(\tau)^{-\frac{1}{100}}$. 
\end{lemma}

\textbf{Proof.} 
By definition of $\delta(\tau)$, we have $|G(0,\tau)| \leq \delta(\tau)$. Moreover, applying Proposition \ref{barrier.rescaled.picture} with $\bar{\tau}=\tau$, we obtain $G_\xi(\xi,\tau)^2 \leq C \, \delta(\tau)$ whenever $G(\xi,\tau) \geq -\frac{1}{\sqrt{2}}$. Putting these facts together, the assertion follows. \\

\begin{lemma}
\label{C2.bound.for.G}
We have $|G_{\xi\xi}(\xi,\tau)| \leq C \, \delta(\tau)^{\frac{1}{8}}$ for $|\xi| \leq \delta(\tau)^{-\frac{1}{100}}$. 
\end{lemma} 

\textbf{Proof.} 
The arguments in the noncompact case carry over unchanged (see \cite{Brendle3}, Lemma 3.8). \\

\begin{lemma}
\label{higher.derivative.bound.for.G}
We have $\big | \frac{\partial^m}{\partial \xi^m} G(\xi,\tau) \big | \leq C(m)$ for $|\xi| \leq \delta(\tau)^{-\frac{1}{100}}$. 
\end{lemma}

\textbf{Proof.} 
The arguments in the noncompact case carry over unchanged (see \cite{Brendle3}, Lemma 3.9). \\

\begin{lemma}
\label{integral.estimates.for.G}
We have 
\[|G_\xi(0,\tau)|^4 \leq C \, \delta(\tau)^{\frac{1}{100}} \int_{\{|\xi| \leq \delta(\tau)^{-\frac{1}{100}}\}} e^{-\frac{\xi^2}{4}} \, |G(\xi,\tau)|^2 \, d\xi\] 
and  
\begin{align*} 
\int_{\{|\xi| \leq \delta(\tau)^{-\frac{1}{100}}\}} e^{-\frac{\xi^2}{4}} \, |G_\xi(\xi,\tau)|^4 \, d\xi 
&\leq C \, \delta(\tau)^{\frac{1}{100}} \int_{\{|\xi| \leq \delta(\tau)^{-\frac{1}{100}}\}} e^{-\frac{\xi^2}{4}} \, |G(\xi,\tau)|^2 \, d\xi \\ 
&+ C \, \exp(-\frac{1}{8} \, \delta(\tau)^{-\frac{1}{50}}). 
\end{align*}
\end{lemma}

\textbf{Proof.} 
The arguments in the noncompact case carry over unchanged (see \cite{Brendle3}, Lemma 3.10). \\

\begin{lemma}
\label{evolution.of.G.2}
We have 
\begin{align*} 
&\int_{\{|\xi| \leq \delta(\tau)^{-\frac{1}{100}}\}} e^{-\frac{\xi^2}{4}} \, \Big | G_\tau(\xi,\tau) - G_{\xi\xi}(\xi,\tau) + \frac{1}{2} \, \xi \, G_\xi(\xi,\tau) - G(\xi,\tau) \Big |^2 \, d\xi \\ 
&\leq C \, \delta(\tau)^{\frac{1}{100}} \int_{\{|\xi| \leq \delta(\tau)^{-\frac{1}{100}}\}} e^{-\frac{\xi^2}{4}} \, |G(\xi,\tau)|^2 \, d\xi + C \, \exp(-\frac{1}{8} \, \delta(\tau)^{-\frac{1}{50}}). 
\end{align*} 
\end{lemma}

\textbf{Proof.} 
The arguments in the noncompact case carry over unchanged (see \cite{Brendle3}, Lemma 3.11). \\

We now perform a spectral decomposition for the operator 
\[\mathcal{L} G := G_{\xi\xi} - \frac{1}{2} \, \xi \, G_\xi + G.\] 
This operator is symmetric with respect to the inner product $\|G\|_{\mathcal{H}}^2 = \int_{\mathbb{R}} e^{-\frac{\xi^2}{4}} \, G^2 \, d\xi$. The eigenvalues of this operator are $1-\frac{n}{2}$, and the corresponding eigenfunctions are $H_n(\frac{\xi}{2})$, where $H_n$ is the $n$-th Hermite polynomial. Let us write $\mathcal{H} = \mathcal{H}_+ \oplus \mathcal{H}_0 \oplus \mathcal{H}_-$, where the subspace $\mathcal{H}_+$ is defined as the span of $H_0(\frac{\xi}{2})$ and $H_1(\frac{\xi}{2})$,  the subspace $\mathcal{H}_0$ is defined as the span of $H_2(\frac{\xi}{2})$, and $\mathcal{H}_-$ is the orthogonal complement of $\mathcal{H}_+ \oplus \mathcal{H}_0$. Moreover, let $P_+$, $P_0$, and $P_-$ denote the orthogonal projections associated with the direct sum $\mathcal{H} = \mathcal{H}_+ \oplus \mathcal{H}_0 \oplus \mathcal{H}_-$. The eigenvalues of the operator $-\mathcal{L}$ on $\mathcal{H}_+$ are bounded from above by $-\frac{1}{2}$. Similarly, the eigenvalues of the operator $-\mathcal{L}$ on $\mathcal{H}_-$ are bounded from below by $\frac{1}{2}$.

Let $\chi$ denote a cutoff function satisfying $\chi(s)=1$ for $s \in [-\frac{1}{2},\frac{1}{2}]$, $\chi(s)=0$ for $s \in \mathbb{R} \setminus [-1,1]$, and $s \chi'(s) \leq 0$ for all $s \in \mathbb{R}$. We define
\begin{align*} 
\gamma(\tau) &:= \int_{\mathbb{R}} e^{-\frac{\xi^2}{4}} \, |G(\xi,\tau) \, \chi(\delta(\tau)^{\frac{1}{100}} \xi)|^2 \, d\xi, \\ 
\gamma^+(\tau) &:= \int_{\mathbb{R}} e^{-\frac{\xi^2}{4}} \, |P_+ (G(\xi,\tau) \, \chi(\delta(\tau)^{\frac{1}{100}} \xi))|^2 \, d\xi, \\ 
\gamma^0(\tau) &:= \int_{\mathbb{R}} e^{-\frac{\xi^2}{4}} \, |P_0 (G(\xi,\tau) \, \chi(\delta(\tau)^{\frac{1}{100}} \xi))|^2 \, d\xi, \\ 
\gamma^-(\tau) &:= \int_{\mathbb{R}} e^{-\frac{\xi^2}{4}} \, |P_- (G(\xi,\tau) \, \chi(\delta(\tau)^{\frac{1}{100}} \xi))|^2 \, d\xi.
\end{align*}
Clearly, $\frac{1}{C} \, \gamma(\tau) \leq \gamma^+(\tau) + \gamma^0(\tau) + \gamma^-(\tau) \leq C \, \gamma(\tau)$. Using Lemma \ref{C1.bound.for.G}, we obtain 
\[\gamma(\tau) \leq C \sup_{\{|\xi| \leq \delta(\tau)^{-\frac{1}{100}}\}} |G(\xi,\tau)|^2 \leq C \, \delta(\tau)^{\frac{1}{4}}.\] 
In particular, $\gamma(\tau) \to 0$ as $\tau \to -\infty$. 

We first analyze the evolution of $\gamma^+(\tau)$, $\gamma^0(\tau)$, and $\gamma^-(\tau)$.

\begin{lemma}
\label{discrete.dynamical.system}
We have
\begin{align*} 
&\gamma^+(\tau-1) \leq e^{-1} \, \gamma^+(\tau) + C \, \delta(\tau)^{\frac{1}{200}} \, \sup_{[\tau-1,\tau]} \gamma(\cdot) + C \, \exp(-\frac{1}{64} \, \delta(\tau)^{-\frac{1}{50}}), \\ 
&|\gamma^0(\tau-1) - \gamma^0(\tau)| \leq C \, \delta(\tau)^{\frac{1}{200}} \, \sup_{[\tau-1,\tau]} \gamma(\cdot) + C \, \exp(-\frac{1}{64} \, \delta(\tau)^{-\frac{1}{50}}), \\ 
&\gamma^-(\tau-1) \geq e \, \gamma^-(\tau) - C \, \delta(\tau)^{\frac{1}{200}} \, \sup_{[\tau-1,\tau]} \gamma(\cdot) - C \, \exp(-\frac{1}{64} \, \delta(\tau)^{-\frac{1}{50}}).
\end{align*} 
\end{lemma} 

\textbf{Proof.} 
This follows from Lemma \ref{evolution.of.G.2}. The proof is analogous to the noncompact case (see \cite{Brendle3}, Lemma 3.12). \\

We next analyze the evolution of the function $\rho_{\text{\rm max}}(\tau)$. We begin with a lemma:

\begin{lemma}
\label{second.derivative.of.G.at.the.maximum.point}
Suppose that $-\tau$ is sufficiently large, and $\xi_*$ be the point in space where the function $\xi \mapsto G(\xi,\tau)$ attains its maximum. Then $0 < -G_{\xi\xi}(\xi_*,\tau) \leq C \, \gamma(\tau)^{\frac{1}{4}}$. 
\end{lemma}

\textbf{Proof.} 
Without loss of generality, we may assume that $\xi_* \geq 0$. By Proposition \ref{neck}, every point with $G(\xi,\tau) \geq -\frac{1}{\sqrt{2}}$ lies on an evolving $\varepsilon$-neck if $-\tau$ is sufficiently large (depending on $\varepsilon$). In particular, $2 \, |G_{\xi\xi\xi}(\xi,\tau)| \leq C$ whenever $G(\xi,\tau) \geq -\frac{1}{\sqrt{2}}$. Consequently, $-G_{\xi\xi}(\xi,\tau) \geq -\frac{1}{2} \, G_{\xi\xi}(\xi_*,\tau) > 0$ for all $\xi \in [\xi_* + \frac{1}{C} \, G_{\xi\xi}(\xi_*,\tau),\xi_*]$. Since $G_\xi(\xi_*,\tau) = 0$, it follows that $G_\xi(\xi,\tau) \geq \frac{1}{2C} \, G_{\xi\xi}(\xi_*,\tau)^2$ for all $\xi \leq \xi_* + G_{\xi\xi}(\xi_*,\tau)$. Since $\xi_* \geq 0$ and $G_{\xi\xi}(\xi_*,\tau)$ is very small, we know that $\xi_* + G_{\xi\xi}(\xi_*,\tau) \geq -1$, and consequently $G_\xi(\xi,\tau) \geq \frac{1}{2C} \, G_{\xi\xi}(\xi_*,\tau)^2$ for all $\xi \leq -1$. We distinguish two cases:

\textit{Case 1:} If $G(-3,\tau) \geq 0$, then $G(\xi,\tau) \geq \frac{1}{2C} \, G_{\xi\xi}(\xi_*,\tau)^2$ for $\xi \in [-2,-1]$.

\textit{Case 2:} If $G(-3,\tau) \leq 0$, then $G(\xi,\tau) \leq -\frac{1}{2C} \, G_{\xi\xi}(\xi_*,\tau)^2$ for $\xi \in [-5,-4]$.

In either case, we conclude that $\int_{-5}^5 |G(\xi,\tau)|^2 \, d\xi \geq \frac{1}{4C^2} \, G_{\xi\xi}(\xi_*,\tau)^4$. Consequently, $G_{\xi\xi}(\xi_*,\tau)^4$ is bounded by a large constant times $\gamma(\tau)$. This completes the proof of Lemma \ref{second.derivative.of.G.at.the.maximum.point}. \\

\begin{lemma}
\label{evolution.of.rho_max}
The function $\rho_{\text{\rm max}}(\tau)$ satisfies 
\[\frac{d}{d\tau} \rho_{\text{\rm max}}(\tau) \geq \rho_{\text{\rm max}}(\tau) - C\rho_{\text{\rm max}}(\tau)^2 - C \, \gamma(\tau)^{\frac{1}{4}}\] 
if $-\tau$ is sufficiently large.
\end{lemma}

\textbf{Proof.} 
We compute 
\[\frac{d}{d\tau} \rho_{\text{\rm max}}(\tau) = \frac{1}{2} \, (\sqrt{2}+\rho_{\text{\rm max}}(\tau)) - (\sqrt{2}+\rho_{\text{\rm max}}(\tau))^{-1} + G_{\xi\xi}(\xi_*,\tau),\] 
where $\xi_*$ is the point where the function $\xi \mapsto G(\xi,\tau)$ attains its maximum. Using Lemma \ref{second.derivative.of.G.at.the.maximum.point}, we obtain 
\[\frac{d}{d\tau} \rho_{\text{\rm max}}(\tau) \geq \rho_{\text{\rm max}}(\tau) - C\rho_{\text{\rm max}}(\tau)^2 - C \, \gamma(\tau)^{\frac{1}{4}}.\] 
This proves the assertion. \\

\begin{lemma}
\label{discrete.dynamical.system.for.rho_max} 
We have 
\[\rho_{\text{\rm max}}(\tau-1) \leq e^{-\frac{1}{2}} \, \rho_{\text{\rm max}}(\tau) + C \, \sup_{[\tau-1,\tau]} \gamma(\cdot)^{\frac{1}{4}}\]
if $-\tau$ is sufficiently large.
\end{lemma}

\textbf{Proof.} 
Lemma \ref{evolution.of.rho_max} implies $\frac{d}{d\tau} \rho_{\text{\rm max}}(\tau) \geq \frac{1}{2} \, \rho_{\text{\rm max}}(\tau) - C \, \gamma(\tau)^{\frac{1}{4}}$ if $-\tau$ is sufficiently large. If we integrate this differential inequality, the assertion follows. \\

\begin{lemma}
\label{combined.estimate}
We have
\begin{align*} 
&\gamma^+(\tau-1) + \rho_{\text{\rm max}}(\tau-1)^8 \\ 
&\leq e^{-1} \, (\gamma^+(\tau) + \rho_{\text{\rm max}}(\tau)^8) + C \, \delta(\tau)^{\frac{1}{200}} \, \sup_{[\tau-1,\tau]} \gamma(\cdot) + C \, \exp(-\frac{1}{64} \, \delta(\tau)^{-\frac{1}{50}}). 
\end{align*}
\end{lemma}

\textbf{Proof.} 
By Lemma \ref{discrete.dynamical.system}, we have  
\[\gamma^+(\tau-1) \leq e^{-1} \, \gamma(\tau)^+ + C \, \delta(\tau)^{\frac{1}{200}} \, \sup_{[\tau-1,\tau]} \gamma(\cdot) + C \, \exp(-\frac{1}{64} \, \delta(\tau)^{-\frac{1}{50}})\] 
Moreover, using Lemma \ref{discrete.dynamical.system.for.rho_max} and Young's inequality, we obtain 
\[\rho_{\text{\rm max}}(\tau-1)^8 \leq e^{-1} \, \rho_{\text{\rm max}}(\tau)^8 + C \, \sup_{[\tau-1,\tau]} \gamma(\cdot)^2.\] 
Adding these inequalities, the assertion follows. \\

We now define 
\begin{align*} 
&\Gamma(\bar{\tau}) := \sup_{\tau \leq \bar{\tau}} (\gamma(\tau) + \rho_{\text{\rm max}}(\tau)^8), \\ 
&\Gamma^+(\bar{\tau}) := \sup_{\tau \leq \bar{\tau}} (\gamma^+(\tau) + \rho_{\text{\rm max}}(\tau)^8), \\ 
&\Gamma^0(\bar{\tau}) := \sup_{\tau \leq \bar{\tau}} \gamma^0(\tau), \\ 
&\Gamma^-(\bar{\tau}) := \sup_{\tau \leq \bar{\tau}} \gamma^-(\tau). 
\end{align*} 
Clearly, $\frac{1}{C} \, \Gamma(\bar{\tau}) \leq \Gamma^+(\bar{\tau}) + \Gamma^0(\bar{\tau}) + \Gamma^-(\bar{\tau}) \leq C \, \Gamma(\bar{\tau})$. It follows from Lemma \ref{C1.bound.for.G} that $\gamma(\tau) \leq C \, \delta(\tau)^{\frac{1}{4}}$. Moreover, $\rho_{\text{\rm max}}(\tau) \leq \delta(\tau)$ by definition of $\delta(\tau)$. Putting these facts together gives $\Gamma(\bar{\tau}) \leq C \, \delta(\bar{\tau})^{\frac{1}{4}}$. In particular, $\Gamma(\bar{\tau}) \to 0$ as $\bar{\tau} \to -\infty$. Using Lemma \ref{combined.estimate} and Lemma \ref{discrete.dynamical.system}, we obtain 
\begin{align*} 
&\Gamma^+(\bar{\tau}-1) \leq e^{-1} \, \Gamma^+(\bar{\tau}) + C \, \delta(\bar{\tau})^{\frac{1}{200}} \, \Gamma(\bar{\tau}) + C \, \exp(-\frac{1}{64} \, \delta(\bar{\tau})^{-\frac{1}{50}}), \\ 
&|\Gamma^0(\bar{\tau}-1) - \Gamma^0(\bar{\tau})| \leq C \, \delta(\bar{\tau})^{\frac{1}{200}} \, \Gamma(\bar{\tau}) + C \, \exp(-\frac{1}{64} \, \delta(\bar{\tau})^{-\frac{1}{50}}), \\ 
&\Gamma^-(\bar{\tau}-1) \geq e \, \Gamma^-(\bar{\tau}) - C \, \delta(\bar{\tau})^{\frac{1}{200}} \, \Gamma(\bar{\tau}) - C \, \exp(-\frac{1}{64} \, \delta(\bar{\tau})^{-\frac{1}{50}}).
\end{align*} 
It follows from standard interpolation inequalities that $|G(0,\tau)| \leq C \, \gamma(\tau)^{\frac{1}{4}}$, hence $\sup_{\tau \leq \bar{\tau}} |G(0,\tau)| \leq C \, \Gamma(\bar{\tau})^{\frac{1}{4}}$. Since $\sup_{\tau \leq \bar{\tau}} \rho_{\text{\rm max}}(\tau) \leq \Gamma(\bar{\tau})^{\frac{1}{8}}$, it follows that 
\[\delta(\bar{\tau}) = \sup_{\tau \leq \bar{\tau}} (|G(0,\tau)|+\rho_{\text{\rm max}}(\tau))  \leq C \, \Gamma(\bar{\tau})^{\frac{1}{8}}.\] 
Consequently, $\exp(-\frac{1}{64} \, \delta(\bar{\tau})^{-\frac{1}{50}}) \leq C \, \delta(\bar{\tau})^9 \leq C \, \delta(\bar{\tau}) \Gamma(\bar{\tau})$. Putting these facts together, we conclude that 
\begin{align*} 
&\Gamma^+(\bar{\tau}-1) \leq e^{-1} \, \Gamma^+(\bar{\tau}) + C \, \delta(\bar{\tau})^{\frac{1}{200}} \, \Gamma(\bar{\tau}), \\ 
&|\Gamma^0(\bar{\tau}-1) - \Gamma^0(\bar{\tau})| \leq C \, \delta(\bar{\tau})^{\frac{1}{200}} \, \Gamma(\bar{\tau}), \\ 
&\Gamma^-(\bar{\tau}-1) \geq e \, \Gamma^-(\bar{\tau}) - C \, \delta(\bar{\tau})^{\frac{1}{200}} \, \Gamma(\bar{\tau}).
\end{align*} 
The following lemma is inspired by a lemma of Merle and Zaag (cf. \cite{Merle-Zaag}, Lemma A.1):

\begin{proposition}
\label{discrete.merle.zaag.lemma}
We either have $\Gamma^0(\bar{\tau}) + \Gamma^-(\bar{\tau}) \leq o(1) \, \Gamma^+(\bar{\tau})$, or $\Gamma^+(\bar{\tau}) + \Gamma^-(\bar{\tau}) \leq o(1) \, \Gamma^0(\bar{\tau})$ as $\bar{\tau} \to -\infty$. 
\end{proposition}

\textbf{Proof.} 
By definition, the function $\Gamma^-(\cdot)$ is monotone increasing. This implies $\Gamma^-(\bar{\tau}) \geq \Gamma^-(\bar{\tau}-1) \geq e \, \Gamma^-(\bar{\tau}) - o(1) \, \Gamma(\bar{\tau})$. Thus, $\Gamma^-(\bar{\tau}) \leq o(1) \, \Gamma(\bar{\tau})$. 

Let $I$ denote the set of all positive real numbers $\alpha$ with the property that the set $\{\bar{\tau}: \Gamma^0(\bar{\tau}) < \alpha \, \Gamma^+(\bar{\tau})\}$ is bounded. Moreover, let $J$ denote the set of all positive real numbers $\alpha$ with the property that the set $\{\bar{\tau}: \Gamma^0(\bar{\tau}) \geq \alpha \, \Gamma^+(\bar{\tau})\}$ is unbounded. Note that $I \subset J$.

We claim that, if $\alpha \in J$, then $e^{\frac{1}{2}} \alpha \in J$ and $e^{-\frac{1}{2}} \alpha \in I$. To see this, suppose that $\alpha \in J$. If we choose $-\bar{\tau}_0$ sufficiently large (depending on $\alpha$), then we obtain 
\[\Gamma^+(\bar{\tau}-1) \leq e^{-1} \, \Gamma^+(\bar{\tau}) + \frac{1}{2(1+\alpha)} \, (e^{-\frac{1}{2}}-e^{-1}) \, (\Gamma^+(\bar{\tau})+\Gamma^0(\bar{\tau}))\] 
and 
\[|\Gamma^0(\bar{\tau}-1)-\Gamma^0(\bar{\tau})| \leq \frac{\alpha}{2(1+\alpha)} \, (1-e^{-\frac{1}{2}}) \, (\Gamma^+(\bar{\tau})+\Gamma^0(\bar{\tau}))\] 
for $\bar{\tau} \leq \bar{\tau}_0$. This gives 
\begin{align*} 
&\Gamma^0(\bar{\tau}-1) - e^{\frac{1}{2}} \alpha \, \Gamma^+(\bar{\tau}-1) \\ 
&\geq \Gamma^0(\bar{\tau}) - e^{-\frac{1}{2}} \alpha \, \Gamma^+(\bar{\tau}) - \frac{\alpha}{1+\alpha} \, (1-e^{-\frac{1}{2}}) \, (\Gamma^+(\bar{\tau})+\Gamma^0(\bar{\tau})) \\ 
&= \Big ( 1-\frac{\alpha}{1+\alpha} \, (1-e^{-\frac{1}{2}}) \Big ) \, (\Gamma^0(\bar{\tau})-\alpha \, \Gamma^+(\bar{\tau})) 
\end{align*} 
for all $\bar{\tau} \leq \bar{\tau}_0$. On the other hand, since $\alpha \in J$, the set $\{\bar{\tau}: \Gamma^0(\bar{\tau}) \geq \alpha \, \Gamma^+(\bar{\tau})\}$ is unbounded. Consequently, there exists some $\bar{\tau}_1 \leq \bar{\tau}_0$ such that $\Gamma^0(\bar{\tau}_1)-\alpha \, \Gamma^+(\bar{\tau}_1) \geq 0$. Proceeding inductively, we obtain $\Gamma^0(\bar{\tau}_1-k)-e^{\frac{1}{2}} \alpha \, \Gamma^+(\bar{\tau}_1-k) \geq 0$ for every positive integer $k$. This implies $e^{\frac{1}{2}} \alpha \in J$. Moreover, we have  
\begin{align*} 
&\inf_{\bar{\tau} \in [\bar{\tau}_1-k-1,\bar{\tau}_1-k]} \big ( \Gamma^0(\bar{\tau}) - e^{-\frac{1}{2}} \alpha \, \Gamma^+(\bar{\tau}) \big ) \\ 
&\geq \Gamma^0(\bar{\tau}_1-k-1) - e^{-\frac{1}{2}} \alpha \, \Gamma^+(\bar{\tau}_1-k) \\ 
&\geq \Gamma^0(\bar{\tau}_1-k) - e^{-\frac{1}{2}} \alpha \, \Gamma^+(\bar{\tau}_1-k) \\ 
&- \frac{\alpha}{2(1+\alpha)} \, (1-e^{-\frac{1}{2}}) \, (\Gamma^+(\bar{\tau}_1-k)+\Gamma^0(\bar{\tau}_1-k)) \\ 
&= \Big ( 1-\frac{\alpha}{2(1+\alpha)} \, (1-e^{-\frac{1}{2}}) \Big ) \, (\Gamma^0(\bar{\tau}_1-k)-\alpha \, \Gamma^+(\bar{\tau}_1-k)) \\ 
&+ \frac{\alpha}{2} \, (1-e^{-\frac{1}{2}}) \, \Gamma^+(\bar{\tau}_1-k) \\ 
&\geq 0 
\end{align*} 
for every positive integer $k$. This implies $e^{-\frac{1}{2}} \alpha \in I$. This completes the proof of the claim. 

Using the claim, we conclude that either $J = \emptyset$ or $I = (0,\infty)$. If $I = (0,\infty)$, it follows that $\Gamma^+(\bar{\tau}) \leq o(1) \, \Gamma^0(\bar{\tau})$ as $\bar{\tau} \to -\infty$. On the other hand, if $J = \emptyset$, then $\Gamma^0(\bar{\tau}) \leq o(1) \, \Gamma^+(\bar{\tau})$ as $\bar{\tau} \to -\infty$. This completes the proof of Lemma \ref{discrete.merle.zaag.lemma}. \\

\section{Ruling out the case when the positive modes dominate}

\label{positive.modes.cannot.dominate}

In this section, we will show that the positive modes cannot dominate. Recall that $r_{\text{\rm max}}(t) \geq \sqrt{-2t}$ for all $t$.

\begin{definition} 
Given $0 < \alpha < 1$, we say that condition $(\star_\alpha)$ holds if $r_{\text{\rm max}}(t) \leq \sqrt{-2t} \, (1+O(-t)^{-\alpha})$. 
\end{definition}

\begin{proposition}
\label{consequence.of.maximum.principle}
Suppose that $(\star_\alpha)$ holds for some $0 < \alpha < 1$. If $-t$ is sufficiently large, then $F_z(z,t)^2 \leq C \, (-t)^{-\frac{\alpha}{1-\alpha}}$ whenever $F(z,t) \geq \sqrt{-t}$. 
\end{proposition} 

\textbf{Proof.} 
Let $\psi_a$ denote the functions in Proposition \ref{properties.of.psi_a}. By assumption, $\frac{r_{\text{\rm max}}(t)}{\sqrt{-2t}} \leq 1+O((-t)^{-\alpha})$. Hence, we can find a constant $K$ with the following properties: 
\begin{itemize}
\item $\frac{r_{\text{\rm max}}(t)^2+2t}{r_{\text{\rm max}}(t)^2} \leq \frac{K^\alpha}{100^{1-\alpha} \, r_{\text{\rm max}}(t)^{2\alpha}}$ for $-t \geq K^{\frac{2}{\alpha}}$. 
\item $\psi_a(s) \geq a^{-2} (s^{-2}-1) + \frac{1}{16} \, a^{-4}$ for $s \in [1-\theta,1+\frac{1}{100} \, a^{-2}]$ and $a \geq K$.
\end{itemize}
In the following, we consider a large number $a \geq K$. Using the inequality $r_{\text{\rm max}}(t) \geq \sqrt{-2t}$, we obtain 
\[\frac{r_{\text{\rm max}}(t)}{\sqrt{-2t+Ka^{\frac{2(1-\alpha)}{\alpha}}}} \geq \frac{r_{\text{\rm max}}(t)}{\sqrt{-(2+K^{-1})t}} \geq \sqrt{\frac{2}{2+K^{-1}}} 
\] 
if $-t \geq K^2 a^{\frac{2(1-\alpha)}{\alpha}}$. Moreover, 
\begin{align*} 
&\frac{-2t+Ka^{\frac{2(1-\alpha)}{\alpha}}}{r_{\text{\rm max}}(t)^2} - 1 + \frac{1}{100} \, a^{-2} \\ 
&= \frac{Ka^{\frac{2(1-\alpha)}{\alpha}}}{r_{\text{\rm max}}(t)^2} - \frac{r_{\text{\rm max}}(t)^2+2t}{r_{\text{\rm max}}(t)^2} + \frac{1}{100} \, a^{-2} \\ 
&\geq \frac{Ka^{\frac{2(1-\alpha)}{\alpha}}}{r_{\text{\rm max}}(t)^2} - \frac{K^\alpha}{100^{1-\alpha} \, r_{\text{\rm max}}(t)^{2\alpha}} + \frac{1}{100} \, a^{-2} \\ 
&\geq 0 
\end{align*} 
for $-t \geq K^2  a^{\frac{2(1-\alpha)}{\alpha}}$. Note that in the last step we have used the elementary inequality $x+y \geq x^\alpha y^{1-\alpha}$ with $x := \frac{Ka^{\frac{2(1-\alpha)}{\alpha}}}{r_{\text{\rm max}}(t)^2}$ and $y := \frac{1}{100} \, a^{-2}$.

Consequently, 
\[\frac{r_{\text{\rm max}}(t)}{\sqrt{-2t+Ka^{\frac{2(1-\alpha)}{\alpha}}}} \leq (1-\frac{1}{100} \, a^{-2})^{-\frac{1}{2}} \leq 1+\frac{1}{100} \, a^{-2}\] 
for $-t \geq K^2  a^{\frac{2(1-\alpha)}{\alpha}}$. In other words, 
\[\frac{F(z,t)}{\sqrt{-2t+Ka^{\frac{2(1-\alpha)}{\alpha}}}} \leq \frac{r_{\text{\rm max}}(t)}{\sqrt{-2t+Ka^{\frac{2(1-\alpha)}{\alpha}}}} \leq 1+\frac{1}{100} \, a^{-2}\] 
for all $-t \geq K^2  a^{\frac{2(1-\alpha)}{\alpha}}$ and all $z$. 

In the following, we fix a real number $a \geq K$, and denote by $\mathcal{I}$ the set of all times $t \leq -K^2 a^{\frac{2(1-\alpha)}{\alpha}}$ with the property that 
\[F_z(z,t)^2 < \psi_a \bigg ( \frac{F(z,t)}{\sqrt{-2t+Ka^{\frac{2(1-\alpha)}{\alpha}}}} \bigg )\] 
whenever $F(z,t) \geq r_* a^{-1} \sqrt{-2t+Ka^{\frac{2(1-\alpha)}{\alpha}}}$. It follows from Proposition \ref{neck} that $t \in \mathcal{I}$ if $-t$ is sufficiently large (depending on $a$). 

We claim that $\mathcal{I} = (-\infty,-K^2 a^{\frac{2(1-\alpha)}{\alpha}}]$. Suppose this is false. Let $t_0$ denote the infimum of the set $(-\infty,-K^2 a^{\frac{2(1-\alpha)}{\alpha}}] \setminus \mathcal{I}$. We can find a point $z_0$ such that $F(z_0,t_0) \geq r_* a^{-1} \sqrt{-2t_0+Ka^{\frac{2(1-\alpha)}{\alpha}}}$ and 
\[F_z(z_0,t_0)^2 = \psi_a \bigg ( \frac{F(z_0,t_0)}{\sqrt{-2t_0+Ka^{\frac{2(1-\alpha)}{\alpha}}}} \bigg ).\] 
Clearly, $F_z(z_0,t_0) \neq 0$ since $\psi_a$ is positive; in other words, the function $z \mapsto F(z,t_0)$ does not attain its maximum at $z_0$.

If $F(z_0,t_0) = r_* a^{-1} \sqrt{-2t_0+Ka^{\frac{2(1-\alpha)}{\alpha}}}$, then 
\[F(z_0,t_0)^2 \leq 1 < \psi_a(r_* a^{-1}) = \psi_a \bigg ( \frac{F(z_0,t_0)}{\sqrt{-2t_0+Ka^{\frac{2(1-\alpha)}{\alpha}}}} \bigg ),\]
which contradicts our choice of $(z_0,t_0)$. Thus, we conclude that $F(z_0,t_0) > r_* a^{-1} \sqrt{-2t_0+Ka^{\frac{2(1-\alpha)}{\alpha}}}$. 

Since $F_z(z_0,t_0) \neq 0$, we can find a smooth function $u(r,t)$ such that $F_z(z,t)^2 = u(F(z,t),t)$ in a neighborhood of the point $(z_0,t_0)$. The function $u$ is defined locally in a neighborhood of $(r_0,t_0)$, where $r_0 := F(z_0,t_0)$. Clearly, 
\[u(r_0,t_0) = \psi_a \bigg ( \frac{r_0}{\sqrt{-2t_0+Ka^{\frac{2(1-\alpha)}{\alpha}}}} \bigg ).\] 
Moreover, since $(-\infty,t_0) \subset \mathcal{I}$, we obtain 
\[u(r,t) < \psi_a \bigg ( \frac{r}{\sqrt{-2t+Ka^{\frac{2(1-\alpha)}{\alpha}}}} \bigg )\] 
for $t < t_0$. As in \cite{Brendle3}, the function $u(r,t)$ satisfies 
\[u_t = uu_{rr} - \frac{1}{2} \, u_r^2 + r^{-2} \, (1-u) \, (ru_r+2u).\] 
On the other hand, the function 
\[\tilde{\Psi}_a(r,t) := \psi_a \bigg ( \frac{r}{\sqrt{-2t+Ka^{\frac{2(1-\alpha)}{\alpha}}}} \bigg )\]
satisfies 
\[\tilde{\Psi}_{a,t} > \tilde{\Psi}_a \tilde{\Psi}_{a,rr} - \frac{1}{2} \, \tilde{\Psi}_{a,r}^2 + r^{-2} \, (1-\tilde{\Psi}_a) \, (r \tilde{\Psi}_{a,r}+2\tilde{\Psi}_a).\] 
This contradicts the parabolic maximum principle. 

To summarize, we have shown that $\mathcal{I} = (-\infty,-K^2 a^{\frac{2(1-\alpha)}{\alpha}}]$. Consequently, $F_z(z,t)^2 \leq C a^{-2}$ whenever $-t \geq K^2  a^{\frac{2(1-\alpha)}{\alpha}}$ and $F(z,t) \geq \sqrt{-t}$. Putting $t = - K^2  a^{\frac{2(1-\alpha)}{\alpha}}$, the assertion follows. This completes the proof of Proposition \ref{consequence.of.maximum.principle}. \\

\begin{corollary}
\label{diameter.bound}
Suppose that $(\star_\alpha)$ holds for some $0 < \alpha < 1$. Then $\text{\rm diam}(S^3,g(t)) \geq \frac{1}{C} \, (-t)^{\frac{1}{2(1-\alpha)}}$ if $-t$ is sufficiently large.
\end{corollary}

\textbf{Proof.} 
We know that the maximum value of $F(z,t)$ is at least $\sqrt{-2t}$. Moreover, Proposition \ref{consequence.of.maximum.principle} implies that $|F_z(z,t)| \leq C \, (-t)^{-\frac{\alpha}{2(1-\alpha)}}$ whenever $F(z,t) \geq \sqrt{-t}$. From this, the assertion follows. \\

\begin{proposition}
\label{key.estimate}
Suppose that $(\star_\alpha)$ holds for some $0 < \alpha < 1$. Moreover, suppose that $(q_0,t_0)$ in space-time with the property that the sphere of symmetry through the point $(q_0,t_0)$ has radius at least $\sqrt{-2t_0}$. If $-t_0$ is sufficiently large, then $-(F^2)_{zz} \leq C \, (-t_0)^{-(1+\frac{\alpha^2}{200})\alpha}$ at the point $(q_0,t_0)$.
\end{proposition}

\textbf{Proof.} We denote by $\tilde{F}(z,t)$ the radius of the sphere of symmetry which has signed distance $z$ from the point $q_0$. The function $\tilde{F}$ satisfies the evolution equation 
\begin{align*} 
\tilde{F}_t(z,t) 
&= \tilde{F}_{zz}(z,t) - \tilde{F}(z,t)^{-1} \, (1+\tilde{F}_z(z,t)^2) \\ 
&+ 2 \, \tilde{F}_z(z,t) \, \bigg [ \tilde{F}(0,t)^{-1} \, \tilde{F}_z(0,t) - \int_0^z \frac{\tilde{F}_z(z',t)^2}{\tilde{F}(z',t)^2} \, dz' \bigg ]. 
\end{align*} 
In particular, 
\[\tilde{F}_t(0,t) = \tilde{F}_{zz}(0,t) - \tilde{F}(0,t)^{-1} \, (1-\tilde{F}_z(0,t)^2).\] 
Using the inequality $\tilde{F}_{zz}(0,t) \leq 0$, we obtain 
\[\tilde{F}_t(0,t) \leq -\tilde{F}(0,t)^{-1} \, (1-\tilde{F}_z(0,t)^2),\] 
hence 
\[-\frac{1}{2} \, \frac{d}{dt} (\tilde{F}(0,t)^2) \geq 1-\tilde{F}_z(0,t)^2.\] 
It follows from Proposition \ref{consequence.of.maximum.principle} that $\tilde{F}_z(0,t)^2 \leq C \, (-t)^{-\frac{\alpha}{1-\alpha}} \leq C \, (-t)^{-\alpha}$ whenever $\tilde{F}(0,t)^2 \geq -t$. This implies 
\[-\frac{1}{2} \, \frac{d}{dt} (\tilde{F}(0,t)^2) \geq 1-C \, (-t)^{-\alpha} > \frac{1}{2}\] 
whenever $\tilde{F}(0,t)^2 \geq -t$ and $-t$ is sufficiently large. By assumption, $\tilde{F}(0,t_0)^2 \geq -2t_0$. Hence, if $-t_0$ is sufficiently large, then a standard continuity argument gives 
\[\tilde{F}(0,t)^2 \geq (-2t) \, (1 - C \, (-t)^{-\alpha})\] 
for all $t \leq t_0$. In the following, we put $\varepsilon := \frac{\alpha^2}{100}$. Using the estimate for $\tilde{F}_z(z,t)^2$  in Proposition \ref{consequence.of.maximum.principle}, we obtain 
\[\tilde{F}(z,t)^2 \geq (-2t) \, (1 - C \, (-t)^{-\frac{\alpha}{8}})\] 
for all $t \leq t_0$ and all $z \in [-(-t)^{\frac{1+\varepsilon}{2}},(-t)^{\frac{1+\varepsilon}{2}}]$. On the other hand, the condition $(\star_\alpha)$ gives 
\[\tilde{F}(z,t)^2 \leq (-2t) \, (1 + C \, (-t)^{-\alpha})\] 
for all $t \leq t_0$ and all $z$. 

We next consider the parabolic cylinder 
\[Q := [-(-t_0)^{\frac{1+\varepsilon}{2}},(-t_0)^{\frac{1+\varepsilon}{2}}] \times [t_0-(-t_0)^{1+\varepsilon},t_0].\] 
Moreover, we define a function $\tilde{H}$ by 
\[\tilde{H}(z,t) := \frac{1}{2} \, \tilde{F}(z,t)^2+t.\] 
By assumption, $\tilde{H}(0,t_0) \geq 0$. Moreover, the preceding arguments imply that we can find a positive constant $L$ such that 
\[-L \, (-t)^{1-\frac{\alpha}{8}} \leq \tilde{H}(z,t) \leq L \, (-t)^{1-\alpha}\] 
in $Q$. In particular, 
\[-L \, (-2t_0)^{(1+\varepsilon)(1-\frac{\alpha}{8})} \leq \tilde{H}(z,t) \leq L \, (-2t_0)^{(1+\varepsilon)(1-\alpha)}\] 
in $Q$. 

The function $\tilde{H}$ satisfies an equation of the form
\[\tilde{H}_t(z,t) - \tilde{H}_{zz}(z,t) = -S(z,t),\] 
where the source term $S$ is defined by 
\[S(z,t) := 2 \, \tilde{F}_z(z,t)^2 - 2 \, \tilde{F}(z,t) \, \tilde{F}_z(z,t) \, \bigg [ \tilde{F}(0,t)^{-1} \, \tilde{F}_z(0,t) - \int_0^z \frac{\tilde{F}_z(z',t)^2}{\tilde{F}(z',t)^2} \, dz' \bigg ].\] 
By Proposition \ref{consequence.of.maximum.principle}, we have $\tilde{F}_z(z,t)^2 \leq C \, (-t)^{-\frac{\alpha}{1-\alpha}}$ at each point in $Q$. This implies 
\[|S(z,t)| \leq C \, (-t)^{-\frac{\alpha}{1-\alpha}} \leq C \, (-t_0)^{-\frac{\alpha}{1-\alpha}}\] 
at each point in $Q$. Moreover, the higher derivatives of $\tilde{F}$ satisfy the estimate $|\frac{\partial^m}{\partial z^m} \tilde{F}(z,t)| \leq C(m) \, (-t_0)^{-\frac{m-1}{2}}$ in the parabolic cylinder $[-(-t_0)^{\frac{1}{2}},(-t_0)^{\frac{1}{2}}] \times [2t_0,t_0]$. (This follows from the pointwise curvature derivative estimate.) This implies $|\frac{\partial^m}{\partial z^m} S(z,t)| \leq C(m) \, (-t_0)^{-\frac{m}{2}}$ in the parabolic cylinder $[-(-t_0)^{\frac{1}{2}},(-t_0)^{\frac{1}{2}}] \times [2t_0,t_0]$. Using standard interpolation inequalities, we obtain 
\begin{align*} 
&|\frac{\partial}{\partial z} S(z,t)| \leq C \, (-t_0)^{-\frac{1}{2}-\frac{\alpha}{1-\alpha}+\varepsilon}, \\ 
&|\frac{\partial^2}{\partial z^2} S(z,t)| \leq C \, (-t_0)^{-1-\frac{\alpha}{1-\alpha}+\varepsilon}, \\ 
&|\frac{\partial}{\partial t} S(z,t)| \leq C \, (-t_0)^{-1-\frac{\alpha}{1-\alpha}+\varepsilon} 
\end{align*} 
in the parabolic cylinder $[-(-t_0)^{\frac{1}{2}},(-t_0)^{\frac{1}{2}}] \times [2t_0,t_0]$.

We now introduce two auxiliary functions $\tilde{H}^{(1)}$ and $\tilde{H}^{(2)}$ on the parabolic cylinder $Q$. Let $\tilde{H}^{(1)}$ denote the solution of the linear heat equation 
\[\tilde{H}^{(1)}_t(z,t) - \tilde{H}^{(1)}_{zz}(z,t) = S(z,t)\] 
on $Q$ with Dirichlet boundary condition $\tilde{H}^{(1)} = 0$ on the parabolic boundary of $Q$. Moreover, let $\tilde{H}^{(2)}$ denote the solution of the linear heat equation 
\[\tilde{H}^{(2)}_t(z,t) - \tilde{H}^{(2)}_{zz}(z,t) = 0\] 
on $Q$ with Dirichlet boundary condition $\tilde{H}^{(2)} = L \, (-2t_0)^{(1+\varepsilon)(1-\alpha)} - \tilde{H}$ on the parabolic boundary of $Q$.

Clearly, $\tilde{H}^{(2)}$ is nonnegative, and we have 
\[\tilde{H}^{(1)}(z,t) + \tilde{H}^{(2)}(z,t) + \tilde{H}(z,t) = L \, (-2t_0)^{(1+\varepsilon)(1-\alpha)}\] 
at each point in $Q$. In particular, 
\[\tilde{H}^{(1)}_{zz}(0,t_0) + \tilde{H}^{(2)}_{zz}(0,t_0) + \tilde{H}_{zz}(0,t_0) = 0.\] 
Therefore, in order to estimate $|\tilde{H}_{zz}(0,t_0)|$, it suffices to bound $|\tilde{H}^{(1)}_{zz}(0,t_0)|$ and $|\tilde{H}^{(2)}_{zz}(0,t_0)|$. 

We begin with the term $|\tilde{H}^{(1)}_{zz}(0,t_0)|$. Using the estimate $|S(z,t)| \leq C \, (-t_0)^{-\frac{\alpha}{1-\alpha}}$ and the maximum principle, we obtain 
\[|\tilde{H}^{(1)}(z,t)| \leq C \, (-t_0)^{1-\frac{\alpha}{1-\alpha}+\varepsilon}\] 
in $Q$. Using standard interior estimates for parabolic equations in the parabolic cylinder $[-(-t_0)^{\frac{1}{2}},(-t_0)^{\frac{1}{2}}] \times [2t_0,t_0]$, we conclude that 
\begin{align*} 
&|\tilde{H}^{(1)}_{zz}(0,t_0)| \\ 
&\leq C \, \sup_{[(-t_0)^{\frac{1}{2}},(-t_0)^{\frac{1}{2}}] \times [2t_0,t_0]} \Big ( (-t_0)^{-1} \, |\tilde{H}^{(1)}| + |S| + (-t_0)^{\frac{1}{2}} \, |\frac{\partial}{\partial z} S| + (-t_0) \, |\frac{\partial}{\partial t} S| \Big ) \\ 
&\leq C \, (-t_0)^{-\frac{\alpha}{1-\alpha}+\varepsilon}. 
\end{align*} 
In the next step, we estimate the term $|\tilde{H}^{(2)}_{zz}(0,t_0)|$. Using the inequality $-\tilde{H}(0,t_0) \leq 0$ together with the estimate 
\[-\tilde{H}^{(1)}(0,t_0) \leq C \, (-t_0)^{1-\frac{\alpha}{1-\alpha}+\varepsilon} \leq C \, (-t_0)^{(1+\varepsilon)(1-\alpha)},\] 
we obtain 
\begin{align*} 
\tilde{H}^{(2)}(0,t_0) 
&= L \, (-2t_0)^{(1+\varepsilon)(1-\alpha)} - \tilde{H}^{(1)}(0,t_0)-\tilde{H}(0,t_0) \\ 
&\leq C \, (-t_0)^{(1+\varepsilon)(1-\alpha)}. 
\end{align*} 
Moreover, we have $\tilde{H}^{(2)} = L \, (-2t_0)^{(1+\varepsilon)(1-\alpha)} - \tilde{H} \leq C \, (-t_0)^{(1+\varepsilon)(1-\frac{\alpha}{8})}$ on the parabolic boundary of $Q$. Hence, applying Proposition \ref{1d.heat.equation} gives 
\begin{align*} 
&(-t_0)^{1+\varepsilon} \, |\tilde{H}^{(2)}_{zz}(0,0)| \\ 
&\leq C \mu^{-2} \, \tilde{H}^{(2)}(0,t_0) + C e^{-\frac{1}{8\mu}} \, \sup_{\{-(-t_0)^{\frac{1+\varepsilon}{2}},(-t_0)^{\frac{1+\varepsilon}{2}}\} \times [t_0-(-t_0)^{1+\varepsilon},t_0]} \tilde{H}^{(2)} \\ 
&\leq C \mu^{-2} \, (-t_0)^{(1+\varepsilon)(1-\alpha)} + C e^{-\frac{1}{8\mu}} (-t_0)^{(1+\varepsilon)(1-\frac{\alpha}{8})}, 
\end{align*}
where $\mu \in (0,1)$ can be chosen arbitrarily. Putting $\mu := (-t_0)^{-\varepsilon^2}$ yields 
\[|\tilde{H}^{(2)}_{zz}(0,0)| \leq C \, (-t_0)^{-(1+\varepsilon)\alpha+2\varepsilon^2}.\] 
Putting these facts together, we conclude that 
\begin{align*} 
|\tilde{H}_{zz}(0,t_0)| 
&\leq |\tilde{H}^{(1)}_{zz}(0,t_0)| + |\tilde{H}^{(2)}_{zz}(0,t_0)| \\ 
&\leq C \, (-t_0)^{-\frac{\alpha}{1-\alpha}+\varepsilon} + C \, (-t_0)^{-(1+\varepsilon)\alpha+2\varepsilon^2} \\ 
&\leq C \, (-t_0)^{-(1+\frac{\varepsilon}{2})\alpha}. 
\end{align*}
This gives $-(\tilde{F}^2)_{zz} \leq C \, (-t_0)^{-(1+\frac{\varepsilon}{2})\alpha}$ at the point $(0,t_0)$. This completes the proof of Proposition \ref{key.estimate}. \\

\begin{corollary}
\label{improved.decay.1}
Suppose that $(\star_\alpha)$ holds for some $0 < \alpha < 1$. If $-t$ is sufficiently large, then $-\frac{1}{2} \, \frac{d}{dt} (r_{\text{\rm max}}(t)^2) \leq 1 + C \, (-t)^{-(1+\frac{\alpha^2}{200})\alpha}$.
\end{corollary}

\textbf{Proof.} Consider the point where the radius is maximal. At that point, $F = r_{\text{\rm max}}(t) \geq \sqrt{-2t}$, $F_z = 0$, and $-(F^2)_{zz} \leq C \, (-t)^{-(1+\frac{\alpha^2}{200})\alpha}$ by Proposition \ref{key.estimate}. Using the evolution equation for $F$, we obtain $-\frac{1}{2} \, \frac{\partial}{\partial t} (F^2) \leq 1 + C \, (-t)^{-(1+\frac{\alpha^2}{200})\alpha}$ at the point where the radius is maximal. From this, the assertion follows. \\

\begin{corollary}
\label{improved.decay.2}
Suppose that $(\star_\alpha)$ holds for some $0 < \alpha < 1$. If $0 < \tilde{\alpha} < \min \{(1+\frac{\alpha^2}{200})\alpha,1\}$, then $(\star_{\tilde{\alpha}})$ holds. 
\end{corollary}

\textbf{Proof.} By Corollary \ref{improved.decay.1}, we have $-\frac{1}{2} \, \frac{d}{dt} (r_{\text{\rm max}}(t)^2) \leq 1 + C \, (-t)^{-\tilde{\alpha}}$. Integrating this differential inequality gives $r_{\text{\rm max}}(t)^2 \leq -2t + C \, (-t)^{1-\tilde{\alpha}}$. Consequently, $r_{\text{\rm max}}(t) \leq \sqrt{-2t} \, (1+C \, (-t)^{-\tilde{\alpha}})$. Thus, $(\star_{\tilde{\alpha}})$ holds. This completes the proof of Corollary \ref{improved.decay.2}. \\

After these preparations, we can now rule out the case that the positive modes dominate. Suppose that $\Gamma^0(\bar{\tau}) + \Gamma^-(\bar{\tau}) \leq o(1) \, \Gamma^+(\bar{\tau})$. The results in Section \ref{asymptotics.near.cylinder} imply $\Gamma^+(\bar{\tau}-1) \leq e^{-1} \, \Gamma^+(\bar{\tau}) + C \, \delta(\bar{\tau})^{\frac{1}{200}} \, \Gamma^+(\bar{\tau})$. Iterating this inequality gives $\Gamma^+(\bar{\tau}) \leq O(e^{\frac{\bar{\tau}}{2}})$. In particular, $\rho_{\text{\rm max}}(\tau) \leq O(e^{\frac{\tau}{16}})$. Equivalently, $r_{\text{\rm max}}(t) \leq \sqrt{-2t} \, (1+O((-t)^{-\frac{1}{16}})$. Therefore, $(\star_\alpha)$ holds for $\alpha=\frac{1}{16}$. Iterating Corollary \ref{improved.decay.2} finitely many times, we conclude that $(\star_\alpha)$ holds for each $0 < \alpha < 1$. Using Corollary \ref{diameter.bound}, we obtain 
\[\liminf_{t \to -\infty} (-t)^{-\frac{1}{2(1-\alpha)}} \, \text{\rm diam}(S^3,g(t)) > 0\] 
for each $0 < \alpha < 1$. On the other hand, standard estimates for the change of distances under Ricci flow imply $-\frac{d}{dt} \text{\rm diam}(S^3,g(t)) \leq C \, \sqrt{R_{\text{\rm max}}(t)}$. Since $R_{\text{\rm max}}(t)$ is uniformly bounded from above by Hamilton's Harnack inequality, we conclude that 
\[\limsup_{t \to -\infty} (-t)^{-1} \, \text{\rm diam}(S^3,g(t)) < \infty.\] 
This is a contradiction. Thus, we have ruled out the case $\Gamma^0(\bar{\tau}) + \Gamma^-(\bar{\tau}) \leq o(1) \, \Gamma^+(\bar{\tau})$. 

\section{Analysis of the case when the neutral mode dominates}

\label{neutral.mode.dominates}

In view of the preceding discussion, we now focus on the case $\Gamma^+(\bar{\tau}) + \Gamma^-(\bar{\tau}) \leq o(1) \, \Gamma^0(\bar{\tau})$. Recall that $\|G\|_{\mathcal{H}}^2 := \int_{\mathbb{R}} e^{-\frac{\xi^2}{4}} \, G^2 \, d\xi$. Let us define $\|G\|_{\mathcal{D}}^2 := \int_{\mathbb{R}} e^{-\frac{\xi^2}{4}} \, (G^2+G_\xi^2) \, d\xi$. Note that $\|G\|_{\mathcal{D}}^2 = \langle G,(2-\mathcal{L})G \rangle_{\mathcal{H}}$ if $G$ is compactly supported.

Recall that the subspace $\mathcal{H}_0$ is one-dimensional and is spanned by the second Hermite polynomial $H_2(\frac{\xi}{2}) = \xi^2-2$. We consider the projection of the function 
\[\hat{G}(\xi,\tau) := G(\xi,\tau) \, \chi(\delta(\tau)^{\frac{1}{100}} \xi)\] 
to the subspace $\mathcal{H}_0$. More precisely, we write 
\[P_0(\hat{G}(\xi,\tau)) = \sqrt{2} \, \alpha(\tau) \, (\xi^2-2),\] 
where 
\[\alpha(\tau) := \frac{1}{16\sqrt{2\pi}} \int_{\mathbb{R}} e^{-\frac{\xi^2}{4}} \, (\xi^2-2) \, \hat{G}(\xi,\tau) \, d\xi.\] 
Furthermore, we define 
\[A(\bar{\tau}) := \sup_{\tau \leq \bar{\tau}} |\alpha(\tau)|.\] 
Clearly, $\frac{1}{C} \, A(\bar{\tau})^2 \leq \Gamma^0(\bar{\tau}) \leq C \, A(\bar{\tau})^2$ for some constant $C$. This implies $\frac{1}{C} \, A(\bar{\tau})^2 \leq \Gamma(\bar{\tau}) \leq C \, A(\bar{\tau})^2$. Since $\delta(\bar{\tau}) \leq C \, \Gamma(\bar{\tau})^{\frac{1}{8}}$, we conclude that $\delta(\bar{\tau}) \leq C \, A(\bar{\tau})^{\frac{1}{4}}$.

\begin{lemma}
\label{neutral.mode.dominates.positive.mode}
We have $\|P_+ \hat{G}(\xi,\tau)\|_{\mathcal{H}} \leq o(1) \, A(\tau)$.
\end{lemma}

\textbf{Proof.} 
We have 
\[\|P_+ \hat{G}(\xi,\tau)\|_{\mathcal{H}}^2 \leq \Gamma^+(\tau) \leq o(1) \, \Gamma(\tau) \leq o(1) \, A(\tau)^2.\] 
This proves the assertion. \\

\begin{lemma}
\label{neutral.mode.dominates.negative.mode}
We have $\|P_- \hat{G}(\xi,\tau)\|_{\mathcal{H}} \leq C \, \delta(\tau)^{\frac{1}{400}} \, A(\tau)$.
\end{lemma}

\textbf{Proof.}
Recall that $\Gamma^-(\tau-1) \geq e \, \Gamma^-(\tau) - C \, \delta(\tau)^{\frac{1}{200}} \, \Gamma(\tau)$. Since $\Gamma^-(\cdot)$ is monotone increasing, we obtain $\Gamma^-(\tau) \leq C \, \delta(\tau)^{\frac{1}{200}} \, \Gamma(\tau)$. Consequently, 
\[\|P_- \hat{G}(\xi,\tau)\|_{\mathcal{H}}^2 \leq \Gamma^-(\tau) \leq C \, \delta(\tau)^{\frac{1}{200}} \, \Gamma(\tau) \leq C \, \delta(\tau)^{\frac{1}{200}} \, A(\tau)^2.\] 
This proves the assertion. \\

\begin{lemma}
\label{neutral.mode.dominates.negative.mode.2}
We have $\|P_- \hat{G}(\xi,\tau)\|_{\mathcal{D}} \leq C \, \delta(\tau)^{\frac{1}{400}} \, A(\tau)$.
\end{lemma}

\textbf{Proof.} 
Lemma \ref{evolution.of.G.2} implies $\frac{\partial}{\partial \tau} G(\xi,\tau) = \mathcal{L} G(\xi,\tau) + E(\xi,\tau)$, where the source term $E(\xi,\tau)$ satisfies 
\[\int_{\{|\xi| \leq \delta(\tau)^{-\frac{1}{100}}\}} e^{-\frac{\xi^2}{4}} \, E(\xi,\tau)^2 \, d\xi \leq C \, \delta(\tau)^{\frac{1}{100}} \, A(\tau)^2.\] 
Using Lemma \ref{derivative.of.delta}, we obtain $\frac{\partial}{\partial \tau} \hat{G}(\xi,\tau) = \mathcal{L} \hat{G}(\xi,\tau) + \hat{E}(\xi,\tau)$, where the source term $\hat{E}(\xi,\tau)$ satisfies 
\[\int_{\mathbb{R}} e^{-\frac{\xi^2}{4}} \, \hat{E}(\xi,\tau)^2 \, d\xi \leq C \, \delta(\tau)^{\frac{1}{100}} \, A(\tau)^2.\]  
Duhamel's principle gives 
\[\hat{G}(\cdot,\tau) = e^{\mathcal{L}} \, \hat{G}(\cdot,\tau-1) + \int_{\tau-1}^\tau e^{(\tau'-\tau+1)\mathcal{L}} \, \hat{E}(\cdot,\tau') \, d\tau',\] 
where the exponential refers to the semigroup of operators generated by $\mathcal{L}$. This implies 
\[P_- \hat{G}(\cdot,\tau) = e^{\mathcal{L}} \,  P_- \hat{G}(\cdot,\tau-1) + \int_{\tau-1}^\tau e^{(\tau'-\tau+1)\mathcal{L}} \, P_- \hat{E}(\cdot,\tau') \, d\tau'.\] 
Using Lemma \ref{neutral.mode.dominates.negative.mode} together with the estimate for $\|\hat{E}(\cdot,\tau)\|_{\mathcal{H}}$, we obtain   
\begin{align*} 
&\|P_- \hat{G}(\cdot,\tau)\|_{\mathcal{D}} \\ 
&\leq C \, \|P_- \hat{G}(\cdot,\tau-1)\|_{\mathcal{H}} + C \int_{\tau-1}^\tau (\tau'-\tau+1)^{-\frac{1}{2}} \, \|\hat{E}(\cdot,\tau')\|_{\mathcal{H}} \, d\tau' \\ 
&\leq C \, \delta(\tau-1)^{\frac{1}{400}} \, A(\tau-1) + C \int_{\tau-1}^\tau (\tau'-\tau+1)^{-\frac{1}{2}} \, \delta(\tau')^{\frac{1}{200}} \, A(\tau') \, d\tau' \\ 
&\leq C \, \delta(\tau)^{\frac{1}{400}} \, A(\tau), 
\end{align*} 
where in the last step we have used the fact that the functions $\delta(\cdot)$ and $A(\cdot)$ are monotone increasing. \\

\begin{lemma}
\label{bound.for.integrals}
We have 
\[\int_{\mathbb{R}} e^{-\frac{\xi^2}{4}} \, (1+|\xi|)^4 \, |\hat{G}(\xi,\tau) - \sqrt{2} \, \alpha(\tau) \, (\xi^2-2)|^2 \, d\xi \leq o(1) \, A(\tau)^2\] 
and 
\[\int_{\mathbb{R}} e^{-\frac{\xi^2}{4}} \, (1+|\xi|)^4 \, |\hat{G}_\xi(\xi,\tau) - 2\sqrt{2} \, \alpha(\tau) \, \xi|^2 \, d\xi \leq o(1) \, A(\tau)^2.\] 
\end{lemma}

\textbf{Proof.} 
By Lemma \ref{neutral.mode.dominates.positive.mode}, we have 
\[\int_{\mathbb{R}} e^{-\frac{\xi^2}{4}} \, (1+|\xi|)^4 \, |P_+ \hat{G}(\xi,\tau)|^2 \, d\xi \leq o(1) \, A(\tau)^2\] 
and 
\[\int_{\mathbb{R}} e^{-\frac{\xi^2}{4}} \, (1+|\xi|)^4 \, \Big | \frac{\partial}{\partial \xi} P_+ \hat{G}(\xi,\tau) \Big |^2 \, d\xi \leq o(1) \, A(\tau)^2.\] 
Next, Lemma \ref{neutral.mode.dominates.negative.mode.2} implies 
\[\int_{\mathbb{R}} e^{-\frac{\xi^2}{4}} \, |P_- \hat{G}(\xi,\tau)|^2 \, d\xi \leq C \, \delta(\tau)^{\frac{1}{200}} \, A(\tau)^2\] 
and 
\[\int_{\mathbb{R}} e^{-\frac{\xi^2}{4}} \, \Big | \frac{\partial}{\partial \xi} P_- \hat{G}(\xi,\tau) \Big |^2 \, d\xi \leq C \, \delta(\tau)^{\frac{1}{200}} \, A(\tau)^2.\] 
This directly implies 
\[\int_{\{|\xi| \leq \delta(\tau)^{-\frac{1}{1000}}\}} e^{-\frac{\xi^2}{4}} \, (1+|\xi|)^4 \, |P_- \hat{G}(\xi,\tau)|^2 \, d\xi \leq o(1) \, A(\tau)^2\] 
and 
\[\int_{\{|\xi| \leq \delta(\tau)^{-\frac{1}{1000}}\}} e^{-\frac{\xi^2}{4}} \, (1+|\xi|)^4 \, \Big | \frac{\partial}{\partial \xi} P_- \hat{G}(\xi,\tau) \Big |^2 \, d\xi \leq o(1) \, A(\tau)^2.\] 
Consequently, 
\[\int_{\mathbb{R}} e^{-\frac{\xi^2}{4}} \, (1+|\xi|)^4 \, |P_- \hat{G}(\xi,\tau)|^2 \, d\xi \leq o(1) \, A(\tau)^2\] 
and 
\[\int_{\mathbb{R}} e^{-\frac{\xi^2}{4}} \, (1+|\xi|)^4 \, \Big | \frac{\partial}{\partial \xi} P_- \hat{G}(\xi,\tau) \Big |^2 \, d\xi \leq o(1) \, A(\tau)^2.\] 
Putting these facts together, we conclude that 
\[\int_{\mathbb{R}} e^{-\frac{\xi^2}{4}} \, (1+|\xi|)^4 \, |(P_+ + P_-) \hat{G}(\xi,\tau)|^2 \, d\xi \leq o(1) \, A(\tau)^2\] 
and 
\[\int_{\mathbb{R}} e^{-\frac{\xi^2}{4}} \, (1+|\xi|)^4 \, \Big | \frac{\partial}{\partial \xi} (P_+ + P_-) \hat{G}(\xi,\tau) \Big |^2 \, d\xi \leq o(1) \, A(\tau)^2.\] 
This proves the assertion. \\

\begin{lemma}
\label{derivative.of.G.at.0} 
We have $|\hat{G}_\xi(0,\tau)| \leq o(1) \, A(\tau)$.
\end{lemma}

\textbf{Proof.}
Clearly, $\frac{\partial}{\partial \xi} P_0 \hat{G}(\xi,\tau) \big |_{\xi=0} = 0$. Moreover, Lemma \ref{neutral.mode.dominates.positive.mode} implies $\big | \frac{\partial}{\partial \xi} P_+ \hat{G}(\xi,\tau) \big |_{\xi=0} \big | \leq C \, \|P_+ \hat{G}(\cdot,\tau)\|_{\mathcal{H}} \leq o(1) \, A(\tau)$. Finally, using Lemma \ref{neutral.mode.dominates.negative.mode} and standard interpolation inequalities, we obtain $\big | \frac{\partial}{\partial \xi} P_- \hat{G}(\xi,\tau) \big |_{\xi=0} \big | \leq C \, \|P_- \hat{G}(\cdot,\tau)\|_{\mathcal{H}}^{1-\frac{1}{10000}} \leq o(1) \, A(\tau)$. Putting these facts together, the assertion follows. \\

\begin{proposition}
\label{projection.of.source.term.to.neutral.mode}
The function $G$ satisfies $\frac{\partial}{\partial \tau} G(\xi,\tau) = \mathcal{L} G(\xi,\tau) + E(\xi,\tau)$, where 
\[\int_{\{|\xi| \leq \delta(\tau)^{-\frac{1}{100}}\}} e^{-\frac{\xi^2}{4}} \, E(\xi,\tau) \, (\xi^2-2) \, d\xi = -128\sqrt{2\pi} \, \alpha(\tau)^2 + O(A(\tau)^2).\]  
\end{proposition}

\textbf{Proof.} 
The source term $E(\xi,\tau)$ is given by  
\begin{align*} 
E(\xi,\tau) 
&= -G(\xi,\tau) + \frac{1}{2} \, (\sqrt{2}+G(\xi,\tau)) - (\sqrt{2}+G(\xi,\tau))^{-1} \, (1+G_\xi(\xi,\tau)^2) \\ 
&+ 2 \, G_\xi(\xi,\tau) \, \bigg [ (\sqrt{2}+G(0,\tau))^{-1} \, G_\xi(0,\tau) - \int_0^\xi \frac{G_\xi(\xi',\tau)^2}{(\sqrt{2}+G(\xi',\tau))^2} \, d\xi' \bigg ]. 
\end{align*} 
Let us write 
\[E(\xi,\tau) = -\frac{1}{2\sqrt{2}} \, G(\xi,\tau)^2 - \frac{1}{\sqrt{2}} \, G_\xi(\xi,\tau)^2 + E_1(\xi,\tau) + E_2(\xi,\tau) + E_3(\xi,\tau),\] 
where 
\begin{align*} 
E_1(\xi,\tau) &= -G(\xi,\tau) + \frac{1}{2\sqrt{2}} \, G(\xi,\tau)^2 + \frac{1}{2} \, (\sqrt{2}+G(\xi,\tau)) - (\sqrt{2}+G(\xi,\tau))^{-1}, \\ 
E_2(\xi,\tau) &= \Big [ \frac{1}{\sqrt{2}} - (\sqrt{2}+G(\xi,\tau))^{-1} \Big ] \, G_\xi(\xi,\tau)^2, \\ 
E_3(\xi,\tau) &= 2 \, G_\xi(\xi,\tau) \, \bigg [ (\sqrt{2}+G(0,\tau))^{-1} \, G_\xi(0,\tau) - \int_0^\xi \frac{G_\xi(\xi',\tau)^2}{(\sqrt{2}+G(\xi',\tau))^2} \, d\xi' \bigg ].
\end{align*}
Using Lemma \ref{bound.for.integrals}, we obtain 
\begin{align*} 
&\int_{\mathbb{R}} e^{-\frac{\xi^2}{4}} \, \hat{G}(\xi,\tau)^2 \, (\xi^2-2) \, d\xi \\ 
&= 2\alpha(\tau)^2 \int_{\mathbb{R}} e^{-\frac{\xi^2}{4}} \, (\xi^2-2)^3 \, d\xi + o(1) \, A(\tau)^2 = 256\sqrt{\pi} \, \alpha(\tau)^2 + o(1) \, A(\tau)^2 
\end{align*}
and  
\begin{align*} 
&\int_{\mathbb{R}} e^{-\frac{\xi^2}{4}} \, \hat{G}_\xi(\xi,\tau)^2 \, (\xi^2-2) \, d\xi \\ 
&= 8\alpha(\tau)^2 \int_{\mathbb{R}} e^{-\frac{\xi^2}{4}} \, \xi^2 \, (\xi^2-2) \, d\xi  +  o(1) \, A(\tau)^2 = 128\sqrt{\pi} \, \alpha(\tau)^2 + o(1) \, A(\tau)^2. 
\end{align*} 
Finally, we estimate the terms $E_1(\xi,\tau)$, $E_2(\xi,\tau)$, and $E_3(\xi,\tau)$. The term $E_1(\xi,\tau)$ satisfies the pointwise estimate $|E_1(\xi,\tau)| \leq C \, |G(\xi,\tau)|^3 \leq o(1) \, G(\xi,\tau)^2$ for $|\xi| \leq \delta(\tau)^{-\frac{1}{100}}$. Using Lemma \ref{bound.for.integrals}, we obtain 
\begin{align*} 
&\bigg | \int_{\{|\xi| \leq \delta(\tau)^{-\frac{1}{100}}\}} e^{-\frac{\xi^2}{4}} \, E_1(\xi,\tau) \, (\xi^2-2) \, d\xi \bigg | \\ 
&\leq o(1) \int_{\{|\xi| \leq \delta(\tau)^{-\frac{1}{100}}\}} e^{-\frac{\xi^2}{4}} \, G(\xi,\tau)^2 \, (\xi^2+2) \, d\xi \leq o(1) \, A(\tau)^2. 
\end{align*}
The term $E_2(\xi,\tau)$ satisfies the pointwise estimate $|E_2(\xi,\tau)| \leq C \, |G(\xi,\tau)| \, G_\xi(\xi,\tau)^2 \leq o(1) \, G_\xi(\xi,\tau)^2$ for $|\xi| \leq \delta(\tau)^{-\frac{1}{100}}$. Using Lemma \ref{bound.for.integrals}, we obtain 
\begin{align*} 
&\bigg | \int_{\{|\xi| \leq \delta(\tau)^{-\frac{1}{100}}\}} e^{-\frac{\xi^2}{4}} \, E_2(\xi,\tau) \, (\xi^2-2) \, d\xi \bigg | \\ 
&\leq o(1) \int_{\{|\xi| \leq \delta(\tau)^{-\frac{1}{100}}\}} e^{-\frac{\xi^2}{4}} \, G_\xi(\xi,\tau)^2 \, (\xi^2+2) \, d\xi \leq o(1) \, A(\tau)^2. 
\end{align*} 
To estimate the term $E_3(\xi,\tau)$, we observe that the function $\xi \mapsto G_\xi(\xi,\tau)$ is monotone decreasing. Hence, if $\xi'$ lies between $0$ and $\xi$, then $G_\xi(\xi',\tau)$ lies between $G_\xi(0,\tau)$ and $G_\xi(\xi,\tau)$, and consequently $|G_\xi(\xi',\tau)| \leq \max \{|G_\xi(0,\tau)|,|G_\xi(\xi,\tau)|\}$. Using Lemma \ref{derivative.of.G.at.0}, we conclude that the term $E_3(\xi,\tau)$ satisfies the pointwise estimate 
\begin{align*} 
|E_3(\xi,\tau)| 
&\leq 2 \, |G_\xi(\xi,\tau)| \, \Big [ |G_\xi(0,\tau)| + |\xi| \, G_\xi(0,\tau)^2 + |\xi| \, G_\xi(\xi,\tau)^2 \Big ] \\ 
&\leq o(1) \, |G_\xi(\xi,\tau)| \, A(\tau) + o(1) \, |\xi| \, A(\tau)^2 + o(1) \, |\xi| \, G_\xi(\xi,\tau)^2 
\end{align*}
for $|\xi| \leq \delta(\tau)^{-\frac{1}{100}}$. Using Lemma \ref{bound.for.integrals}, we obtain 
\begin{align*} 
&\bigg | \int_{\{|\xi| \leq \delta(\tau)^{-\frac{1}{100}}\}} e^{-\frac{\xi^2}{4}} \, E_3(\xi,\tau) \, (\xi^2-2) \, d\xi \bigg | \\ 
&\leq o(1) \int_{\{|\xi| \leq \delta(\tau)^{-\frac{1}{100}}\}} e^{-\frac{\xi^2}{4}} \, |G_\xi(\xi,\tau)| \, A(\tau) \, (\xi^2+2) \, d\xi \\ 
&+ o(1) \, \int_{\{|\xi| \leq \delta(\tau)^{-\frac{1}{100}}\}} e^{-\frac{\xi^2}{4}} \, A(\tau)^2 \, |\xi| \, (\xi^2+2) \, d\xi \\ 
&+ o(1) \int_{\{|\xi| \leq \delta(\tau)^{-\frac{1}{100}}\}} e^{-\frac{\xi^2}{4}} \, |G_\xi(\xi,\tau)|^2 \, |\xi| \, (\xi^2+2) \, d\xi \\ 
&\leq o(1) \, A(\tau)^2. 
\end{align*} 
Putting these facts together, the assertion follows. \\

\begin{corollary}
\label{ode.for.alpha} 
The function $\alpha(\tau)$ satisfies $\alpha'(\tau) = -8\alpha(\tau)^2 + o(1) \, A(\tau)^2$. 
\end{corollary}

\textbf{Proof.} 
This follows easily from Proposition \ref{projection.of.source.term.to.neutral.mode}. \\

\begin{corollary}
\label{maximum.of.alpha}
If $-\tau$ is sufficiently large, then $A(\tau) = |\alpha(\tau)|$. 
\end{corollary} 

\textbf{Proof.} 
Suppose that $A(\bar{\tau}) > |\alpha(\bar{\tau})|$ for some time $\bar{\tau}$, where $-\bar{\tau}$ is very large. We can find a time $\tau_* \in (-\infty,\bar{\tau})$ such that $|\alpha(\tau_*)| = A(\bar{\tau})$. By continuity, we can find an open interval $I$ such that $\tau_* \in I$, $I \subset (-\infty,\bar{\tau})$, and $|\alpha(\tau)| \geq \frac{1}{2} \, A(\bar{\tau})$ for all $\tau \in I$. Corollary \ref{ode.for.alpha} now implies $\alpha'(\tau) = -8\alpha(\tau)^2 + o(1) \, A(\tau)^2 \leq -A(\bar{\tau})^2$ for almost all $\tau \in I$. Consequently, the function $\alpha(\tau)$ is strictly monotone decreasing on the interval $I$. This contradicts the fact that the function $|\alpha(\tau)|$ attains a local maximum at $\tau_*$. \\

\begin{corollary}
\label{asymptotics.of.alpha} 
The function $\alpha(\tau)$ satisfies $\alpha'(\tau) = -(8+o(1)) \, \alpha(\tau)^2$. In particular, $\alpha(\tau) = \frac{1}{(8+o(1)) \, \tau} < 0$ if $-\tau$ is sufficiently large.
\end{corollary}

\begin{proposition}
\label{convergence.in.C_loc^infty}
We have $(-\tau) \, G(\xi,\tau) \to -\frac{1}{4\sqrt{2}} \, (\xi^2-2)$ in $C_{\text{\rm loc}}^\infty$.
\end{proposition}

\textbf{Proof.} 
Let us fix an arbitrary constant $L$. Lemma \ref{evolution.of.G.2} implies $\frac{\partial}{\partial \tau} G(\xi,\tau) = \mathcal{L} G(\xi,\tau) + E(\xi,\tau)$, where the source term $E(\xi,\tau)$ satisfies 
\[\int_{\{|\xi| \leq \delta(\tau)^{-\frac{1}{100}}\}} e^{-\frac{\xi^2}{4}} \, E(\xi,\tau)^2 \, d\xi \leq C \, \delta(\tau)^{\frac{1}{100}} \, A(\tau)^2.\] 
Note that $A(\tau) \leq C \, (-\tau)^{-1}$ and $\delta(\tau) \leq C \, A(\tau)^{\frac{1}{4}} \leq C \, (-\tau)^{-\frac{1}{4}}$. Therefore, 
\[\int_{\{|\xi| \leq 4L\}} E(\xi,\tau)^2 \, d\xi \leq C(L) \, (-\tau)^{-2-\frac{1}{400}}.\] 
Using standard interpolation inequalities, we obtain 
\[\|E\|_{C^{2m,m}([-2L,2L] \times [\tau-1,\tau])} \leq C(L) \, (-\tau)^{-1-\frac{1}{1000}}\] 
for any given positive integer $m$. Here, $C^{2m,m}$ denotes the space of functions which are $2m$-times continuously differentiable in space and $m$-times continuously differentiable in time.

Recall that $\mathcal{L}(\xi^2-2) = 0$. Hence, the function $G(\xi,\tau) - \frac{1}{4\sqrt{2} \, \tau} \, (\xi^2-2)$ satisfies 
\begin{align*} 
\frac{\partial}{\partial \tau} \Big ( G(\xi,\tau) - \frac{1}{4\sqrt{2} \, \tau} \, (\xi^2-2) \Big ) 
&= \mathcal{L} \Big ( G(\xi,\tau) - \frac{1}{4\sqrt{2} \, \tau} \, (\xi^2-2) \Big ) \\ 
&+ E(\xi,\tau) + \frac{1}{4\sqrt{2} \, \tau^2} \, (\xi^2-2). 
\end{align*}
Moreover, using Lemma \ref{bound.for.integrals} and Corollary \ref{asymptotics.of.alpha}, we obtain 
\[\int_{\{|\xi| \leq 2L\}} \Big | G(\xi,\tau) - \frac{1}{4\sqrt{2} \, \tau} \, (\xi^2-2) \Big |^2 \, d\xi \leq o(1) \, A(\tau)^2 \leq o(1) \, (-\tau)^{-2}.\] 
Using standard interior estimates for parabolic equations, we conclude that 
\[\Big \| G(\xi,\tau) - \frac{1}{4\sqrt{2} \, \tau} \, (\xi^2-2) \Big \|_{C^{2m}([-L,L])} \leq o(1) \, (-\tau)^{-1}\] 
for any given positive integer $m$. This completes the proof. \\

\begin{corollary}
\label{domain.of.definition.of.G}
The domain of definition of the function $\xi \mapsto G(\xi,\tau)$ is an interval of length at most $o(1) \, (-\tau)$.
\end{corollary} 

\textbf{Proof.} 
Recall that the function $\xi \mapsto G(\xi,\tau)$ is concave, and $G(\xi,\tau) \geq -\sqrt{2}$. Hence, the assertion follows from Proposition \ref{convergence.in.C_loc^infty}. \\

\begin{corollary}
\label{max.attained.near.0}
Let $\varepsilon > 0$ be given. If $-\tau$ is sufficiently large, then the function $\xi \mapsto G(\xi,\tau)$ attains its maximum in the interval $(-\varepsilon,\varepsilon)$.
\end{corollary}

\section{Asymptotics in the intermediate region}

\label{intermediate.region.asymptotics}

We next study the asymptotics in the intermediate region where $|z| \geq M\sqrt{-t}$ for some large constant $M$, and $F(z,t) \geq \theta \sqrt{-2t}$ for some small constant $\theta$. The following result is a consequence of our barrier arguments:

\begin{proposition}
\label{consequence.of.barrier.argument}
Fix a small number $\theta \in (0,\frac{1}{2})$ and a large number $M \geq 10$. If $-t$ is sufficiently large (depending on $\theta$ and $M$), then  
\[F_z(z,t)^2 \leq \frac{M^2+C(\theta)}{M^2-2} \, \frac{1}{2 \log(-t)} \, \Big ( \frac{-2t}{F(z,t)^2} - 1 \Big )\] 
whenever $|z| \geq M\sqrt{-t}$ and $F(z,t) \geq \theta \sqrt{-2t}$. 
\end{proposition}

\textbf{Proof.} 
Let us fix a large number $M$. By Proposition \ref{convergence.in.C_loc^infty}, we have 
\[F(z,t) = \sqrt{-2t} \, \Big ( 1 - \frac{M^2-2}{8 \log(-t)} + o \Big ( \frac{1}{\log(-t)} \Big ) \Big )\] 
and 
\[F_z(z,t)^2 = \frac{M^2}{8 \, (\log(-t))^2} + o \Big ( \frac{1}{(\log(-t))^2} \Big )\]
for $|z|=M\sqrt{-t}$. Hence, if $a$ is sufficiently large, then  
\[F_z(z,t)^2 < \psi_a \Big ( \frac{F(z,t)}{\sqrt{-2t}} \Big )\] 
whenever $\log(-t) \geq \frac{M^2+2}{M^2-2} \, \frac{a^2}{2}$ and $|z|=M\sqrt{-t}$. Using the maximum principle, we obtain 
\[F_z(z,t)^2 < \psi_a \Big ( \frac{F(z,t)}{\sqrt{-2t}} \Big )\] 
whenever $\log(-t) \geq \frac{M^2+2}{M^2-2} \, \frac{a^2}{2}$, $|z| \geq M\sqrt{-t}$, and $F(z,t) \geq r_* a^{-1} \sqrt{-2t}$.

We now fix a small number $\theta \in (0,\frac{1}{2})$. Using Definition 2.6 in \cite{Brendle3}, we obtain $\psi_a(s) \leq a^{-2} (s^{-2}-1) + C(\theta) \, a^{-4}$ for $s \in [\theta,1]$. Consequently,  if $a$ is sufficiently large, then we have
\[F_z(z,t)^2 \leq a^{-2} \, \Big ( \frac{-2t}{F(z,t)^2} - 1 \Big ) + C(\theta) \, a^{-4}\] 
whenever $\log(-t) \geq \frac{M^2+2}{M^2-2} \, \frac{a^2}{2}$, $|z| \geq M\sqrt{-t}$, and $F(z,t) \geq \theta \sqrt{-2t}$. We now put $\log(-t) = \frac{M^2+2}{M^2-2} \, \frac{a^2}{2}$. Hence, if $-t$ is sufficiently large, then we have 
\[F_z(z,t)^2 \leq \frac{M^2+2}{M^2-2} \, \frac{1}{2 \log(-t)} \, \Big ( \frac{-2t}{F(z,t)^2} - 1 \Big ) + \frac{C(\theta)}{(\log(-t))^2}\] 
whenever $|z| \geq M\sqrt{-t}$ and $F(z,t) \geq \theta \sqrt{-2t}$. Now, in the region $|z| \geq M \sqrt{-t}$, we have 
\[\frac{M^2-2}{8 \log(-t)} \leq \frac{-2t}{F(z,t)^2} - 1.\] 
Consequently, 
\[F_z(z,t)^2 \leq \frac{M^2+C(\theta)}{M^2-2} \, \frac{1}{2 \log(-t)} \, \Big ( \frac{-2t}{F(z,t)^2} - 1 \Big )\] 
whenever $|z| \geq M\sqrt{-t}$ and $F(z,t) \geq \theta \sqrt{-2t}$. This proves the assertion. \\

\begin{corollary}
\label{time.derivative.of.F}
Fix a small number $\theta \in (0,\frac{1}{2})$. If $-t$ is sufficiently large (depending on $\theta$), then   
\[-(F^2)_t(z,t) \geq 2 - \frac{C(\theta)}{\log(-t)}\] 
whenever $|z| \geq 10\sqrt{-t}$ and $F(z,t) \geq \theta \sqrt{-2t}$. 
\end{corollary}

\textbf{Proof.}
Recall that $F_{zz}(z,t) \leq 0$ at each point in space-time. Moreover, Corollary \ref{max.attained.near.0} implies $F_z(z,t) \geq 0$ for $z \leq -\sqrt{-t}$ and $F_z(z,t) \leq 0$ for $z \geq \sqrt{-t}$. Using the evolution equation for $F$, we obtain 
\begin{align*} 
F_t(z,t) 
&= F_{zz}(z,t) - F(z,t)^{-1} \, (1-F_z(z,t)^2) \\ 
&- 2 \, F_z(z,t) \int_0^z \frac{F_{zz}(z',t)}{F(z',t)} \, dz' \\ 
&\leq -F(z,t)^{-1} \, (1-F_z(z,t)^2) 
\end{align*} 
for $|z| \geq \sqrt{-t}$. Applying Proposition \ref{consequence.of.barrier.argument} with $M=10$ gives $F_z(z,t)^2 \leq \frac{C(\theta)}{\log(-t)}$ whenever $|z| \geq 10\sqrt{-t}$ and $F(z,t) \geq \theta \sqrt{-2t}$. Putting these facts together, we conclude that 
\[F_t(z,t) \leq -F(z,t)^{-1} \, \Big ( 1-\frac{C(\theta)}{\log(-t)} \Big )\] 
whenever $|z| \geq 10\sqrt{-t}$ and $F(z,t) \geq \theta \sqrt{-2t}$. From this, the assertion follows. \\

\begin{proposition} 
\label{lower.bound} 
Fix a small number $\theta \in (0,\frac{1}{2})$ and a large number $M \geq 10$. If $-t$ is sufficiently large (depending on $\theta$ and $M$), then  
\[F(z,t)^2 \geq -2t - \frac{M^2+C(\theta)}{M^2-2} \, \frac{z^2}{2 \log(-t)}\] 
whenever $|z| \geq M\sqrt{-t}$ and $F(z,t) \geq \theta \sqrt{-2t}$. 
\end{proposition}

\textbf{Proof.} 
Proposition \ref{consequence.of.barrier.argument} implies 
\[\Big | \frac{\partial}{\partial z} \sqrt{-2t-F(z,t)^2} \Big | \leq \sqrt{\frac{M^2+C(\theta)}{M^2-2} \, \frac{1}{2 \log(-t)}}\] 
whenever $|z| \geq M\sqrt{-t}$ and $F(z,t) \geq \theta \sqrt{-2t}$. Moreover, by Proposition \ref{convergence.in.C_loc^infty}, we have 
\[\sqrt{-2t-F(z,t)^2} \leq \sqrt{\frac{1}{2\log(-t)}} \, |z|\] 
for $|z|=M\sqrt{-t}$. Consequently, 
\[\sqrt{-2t-F(z,t)^2} \leq \sqrt{\frac{M^2+C(\theta)}{M^2-2} \, \frac{1}{2 \log(-t)}} \, |z|\] 
whenever $|z| \geq M\sqrt{-t}$ and $F(z,t) \geq \theta \sqrt{-2t}$. From this, the assertion follows. \\

\begin{proposition} 
\label{upper.bound}
Fix a small number $\theta \in (0,\frac{1}{2})$ and a large number $M \geq 20$. If $-t$ is sufficiently large (depending on $\theta$ and $M$), then  
\[F(z,t)^2 \leq -2t - \frac{M^2-C(\theta)}{M^2} \, \frac{z^2}{2 \log(-t)}\] 
whenever $|z| \geq M\sqrt{-t}$ and $F(z,t) \geq \theta \sqrt{-2t}$.
\end{proposition}

\textbf{Proof.} 
Let us fix a point $(z_0,t_0)$ such that $|z_0| \geq M\sqrt{-t_0}$ and $F(z_0,t_0) > \theta \sqrt{-2t_0}$. Let $t_* := -\frac{z_0^2}{M^2}$. Clearly, $t_* \leq t_0$, and $|z_0| \geq M \sqrt{-t} \geq 20 \sqrt{-t}$ for all $t \in [t_*,t_0]$. Hence, Corollary \ref{time.derivative.of.F} implies that 
\[-(F^2)_t(z_0,t) \geq 2 - \frac{C(\theta)}{\log(-t)} > 2\theta^2\] 
for all times $t \in [t_*,t_0]$ satisfying $F(z_0,t) \geq \theta \sqrt{-2t}$. Since $F(z_0,t_0) > \theta \sqrt{-2t_0}$, we conclude that $F(z_0,t) > \theta \sqrt{-2t}$ for all times $t \in [t_*,t_0]$, and furthermore 
\[-(F^2)_t(z_0,t) \geq 2 - \frac{C(\theta)}{\log(-t)} \geq 2 - \frac{C(\theta)}{\log(-t_0)}\] 
for all times $t \in [t_*,t_0]$. Integrating this inequality over $[t_*,t_0]$ gives 
\[F(z_0,t_*)^2 - F(z_0,t_0)^2 \geq 2(t_0-t_*) - C(\theta) \, \frac{(-t_*)}{\log(-t_0)}\]
On the other hand, since $|z_0| = M\sqrt{-t_*}$, Proposition \ref{convergence.in.C_loc^infty} implies 
\begin{align*} 
F(z_0,t_*)^2 
&= -2t_* - (M^2-2) \, \frac{(-t_*)}{2 \log(-t_*)} + o \Big ( \frac{(-t_*)}{\log(-t_*)} \Big ) \\ 
&\leq -2t_* - (M^2-4) \, \frac{(-t_*)}{2 \log(-t_*)}. 
\end{align*}
Putting these facts together, we obtain 
\[F(z_0,t_0)^2 \leq -2t_0 - (M^2-4) \, \frac{(-t_*)}{2 \log(-t_*)} + C(\theta) \, \frac{(-t_*)}{\log(-t_0)}.\] 
Using the identity $t_* = -\frac{z_0^2}{M^2}$, we conclude that 
\[F(z_0,t_0)^2 \leq -2t_0 - \frac{M^2-4}{M^2} \, \frac{z_0^2}{2 \log(-t_*)} + \frac{C(\theta)}{M^2} \, \frac{z_0^2}{\log(-t_0)}.\] 
Finally, Corollary \ref{domain.of.definition.of.G} gives $|z_0| \leq \sqrt{-t_0} \, \log(-t_0)$. This implies $-t_* \leq z_0^2 \leq (-t_0) \, (\log(-t_0))^2$, hence $\log (-t_*) \leq \log(-t_0) + 2 \log\log(-t_0) \leq \frac{M^2-4}{M^2-6} \, \log(-t_0)$. Therefore, 
\[F(z_0,t_0)^2 \leq -2t_0 - \frac{M^2-6}{M^2} \, \frac{z_0^2}{2 \log(-t_0)} + \frac{C(\theta)}{M^2} \, \frac{z_0^2}{\log(-t_0)}.\] 
This completes the proof. \\

Combining Proposition \ref{lower.bound} and Proposition \ref{upper.bound} and sending $M \to \infty$, we can draw the following conclusions:

\begin{corollary}
\label{size.of.intermediate.region} 
Fix a small number $\theta \in (0,\frac{1}{2})$. If $-t$ is sufficiently large (depending on $\theta$), then the set $\{z: F(z,t) \geq \theta \sqrt{-2t}\}$ is an interval $[-\bar{z}_1(\theta,t),\bar{z}_2(\theta,t)]$, and 
\begin{align*} 
\bar{z}_1(\theta,t) &= (2+o(1)) \, \sqrt{1-\theta^2} \, \sqrt{(-t) \log(-t)}, \\ 
\bar{z}_2(\theta,t) &= (2+o(1)) \, \sqrt{1-\theta^2} \, \sqrt{(-t) \log(-t)}. 
\end{align*}
\end{corollary}

\begin{corollary}
\label{asymptotics in.intermediate.region} 
Fix a small number $\theta \in (0,\frac{1}{2})$. If $-t$ is sufficiently large (depending on $\theta$), then 
\[F(z,t)^2 = -2t - \frac{z^2}{2\log(-t)} + o(-t)\] 
for $|z| \leq 2\sqrt{1-\theta^2} \sqrt{(-t) \log(-t)}$.
\end{corollary}

\section{Asymptotics in the tip region}

\label{tip.region.asymptotics}

In this final section, we analyze the asymptotics of the solution near each tip. For each $t$, the function $z \mapsto F(z,t)$ is defined on the interval $[-d_{\text{\rm tip},1}(t),d_{\text{\rm tip},2}(t)]$, where $d_{\text{\rm tip},1}(t)$ and $d_{\text{\rm tip},2}(t)$ denote the distance of the reference point $q$ from each tip. We first derive an asymptotic formula for $d_{\text{\rm tip},1}(t)$ and $d_{\text{\rm tip},2}(t)$.

\begin{proposition}
\label{distance.to.tip}
We have 
\[\lim_{t \to -\infty} \frac{d_{\text{\rm tip},1}(t)}{\sqrt{(-t) \, \log(-t)}} = \lim_{t \to -\infty} \frac{d_{\text{\rm tip},2}(t)}{\sqrt{(-t) \, \log(-t)}} = 2.\]
\end{proposition}

\textbf{Proof.} 
For each $\theta \in (0,\frac{1}{4})$, we have $d_{\text{\rm tip},1}(t) \geq \bar{z}_1(\theta,t)$. Using Corollary \ref{size.of.intermediate.region}, we obtain 
\[\liminf_{t \to -\infty} \frac{d_{\text{\rm tip},1}(t)}{\sqrt{(-t) \, \log(-t)}} \geq \liminf_{t \to -\infty} \frac{\bar{z}_1(\theta,t)}{\sqrt{(-t) \, \log(-t)}} = 2\sqrt{1-\theta^2}\] 
for each $\theta \in (0,\frac{1}{4})$. Moreover, since the function $z \mapsto F(z,t)$ is concave, we have $\frac{1}{2} \, (\bar{z}_1(2\theta,t) + d_{\text{\rm tip},1}(t)) \leq \bar{z}_1(\theta,t)$. Using Corollary \ref{size.of.intermediate.region}, we obtain 
\[\limsup_{t \to -\infty} \frac{d_{\text{\rm tip},1}(t)}{\sqrt{(-t) \, \log(-t)}} \leq \limsup_{t \to -\infty} \frac{2\bar{z}_1(\theta,t) - \bar{z}_1(2\theta,t)}{\sqrt{(-t) \, \log(-t)}} = 4\sqrt{1-\theta^2} - 2\sqrt{1-4\theta^2}\] 
for each $\theta \in (0,\frac{1}{4})$. Sending $\theta \to 0$ gives $\lim_{t \to -\infty} \frac{d_{\text{\rm tip},1}(t)}{\sqrt{(-t) \, \log(-t)}} = 2$. An analogous argument gives $\lim_{t \to -\infty} \frac{d_{\text{\rm tip},2}(t)}{\sqrt{(-t) \, \log(-t)}} = 2$. This completes the proof. \\

Finally, we analyze the asymptotic behavior of the scalar curvature at each tip. We first recall a basic fact about the Bryant soliton:

\begin{lemma}
\label{bryant}
Consider the Bryant soliton, normalized so that the scalar curvature at the tip equals $1$. Let $\gamma$ be a geodesic ray emanating from the tip of the Bryant soliton which is parametrized by arclength. Then $\int_0^\infty \text{\rm Ric}(\gamma'(s),\gamma'(s)) \, ds = 1$.
\end{lemma}

\textbf{Proof.} 
See \cite{Brendle3}, Lemma 4.3. \\

We now continue with the analysis of our ancient solution. 

\begin{lemma} 
\label{derivative.of.distance}
We have $-\frac{d}{dt} d_{\text{\rm tip},1}(t) = (1+o(1)) \, R_{\text{\rm tip},1}(t)^{\frac{1}{2}}$ and $-\frac{d}{dt} d_{\text{\rm tip},2}(t) = (1+o(1)) \, R_{\text{\rm tip},2}(t)^{\frac{1}{2}}$.
\end{lemma}

\textbf{Proof.} 
In the following, we assume that a small number $\delta>0$ is given. Let $p_1,p_2 \in S^3$ denote the two tips, so that $d_{\text{\rm tip},1}(t) = d_{g(t)}(p_1,q)$ and $d_{\text{\rm tip},2}(t) = d_{g(t)}(p_2,q)$, where $q$ is our fixed reference point. Note that $p_1$, $p_2$, and $q$ represent fixed points on the manifold. Let us fix a time $t$, and let $\gamma$ denote the minimizing geodesic from the tip $p_1$ to the reference point $q$ with respect to the metric $g(t)$, so that $\gamma(0) = p_1$ and $\gamma(d_{\text{\rm tip},1}(t)) = q$. In view of Lemma \ref{bryant} and Proposition \ref{rescaling.around.tip}, we can find a large constant $A$ (depending on $\delta$) such that $A \geq 8\delta^{-1}$ and 
\[(1-\delta) \, R(p_1,t)^{\frac{1}{2}} \leq \int_0^{A R(p_1,t)^{-\frac{1}{2}}} \text{\rm Ric}(\gamma'(s),\gamma'(s)) \, ds \leq (1+\frac{\delta}{2}) \, R(p_1,t)^{\frac{1}{2}}\] 
if $-t$ is sufficiently large (depending on $\delta$ and $A$). We next observe that $\gamma$ is part of a minimizing geodesic from $p_1$ to $p_2$ of length $d_{\text{\rm tip},1}(t)+d_{\text{\rm tip},2}(t)$. Moreover, Corollary \ref{scalar.curvature.at.tips.2} and Proposition \ref{distance.to.tip} imply that $A R(p_1,t)^{-\frac{1}{2}} < d_{\text{\rm tip},2}(t)$ if $-t$ is sufficiently large (depending on $A$). Hence, we may apply Theorem 17.4(a) in \cite{Hamilton3} with $\sigma = A R(p_1,t)^{-\frac{1}{2}}$ and $L = d_{\text{\rm tip},1}(t) + A R(p_1,t)^{-\frac{1}{2}} < d_{\text{\rm tip},1}(t)+d_{\text{\rm tip},2}(t)$. This gives 
\[0 \leq \int_{AR(p_1,t)^{-\frac{1}{2}}}^{d_{\text{\rm tip},1}(t)} \text{\rm Ric}(\gamma'(s),\gamma'(s)) \, ds \leq 4A^{-1} \, R(p_1,t)^{\frac{1}{2}}.\] 
Putting these facts together, we obtain 
\[(1-\delta) \, R(p_1,t)^{\frac{1}{2}} \leq \int_0^{d_{\text{\rm tip},1}(t)} \text{\rm Ric}(\gamma'(s),\gamma'(s)) \, ds \leq (1+\frac{\delta}{2}+4A^{-1}) \, R(p_1,t)^{\frac{1}{2}}.\] 
if $-t$ is sufficiently large (depending on $\delta$ and $A$). Thus, we conclude that 
\[(1-\delta) \, R_{\text{\rm tip},1}(t)^{\frac{1}{2}} \leq -\frac{d}{dt} d_{\text{\rm tip},1}(t) \leq (1+\frac{\delta}{2}+4A^{-1}) \, R_{\text{\rm tip},1}(t)^{\frac{1}{2}}\] 
if $-t$ is sufficiently large (depending on $\delta$ and $A$). Since $4A^{-1} \leq \frac{\delta}{2}$ and $\delta>0$ can be chosen arbitrarily small, it follows that $-\frac{d}{dt} d_{\text{\rm tip},1}(t) = (1+o(1)) \, R_{\text{\rm tip},1}(t)^{\frac{1}{2}}$. An analogous argument gives $-\frac{d}{dt} d_{\text{\rm tip},2}(t) = (1+o(1)) \, R_{\text{\rm tip},2}(t)^{\frac{1}{2}}$. \\

\begin{proposition}
\label{scalar.curvature.at.tip}
The scalar curvature at each tip satisfies $R_{\text{\rm tip},1}(t) = (1+o(1)) \, \frac{\log (-t)}{(-t)}$ and $R_{\text{\rm tip},2}(t) = (1+o(1)) \, \frac{\log (-t)}{(-t)}$.
\end{proposition}

\textbf{Proof.} 
Let $\varepsilon \in (0,\frac{1}{2})$ be given. By Lemma \ref{derivative.of.distance}, we know that 
\[(1-\varepsilon) \, \sqrt{R_{\text{\rm tip},1}(t)} \leq -\frac{d}{dt} d_{\text{\rm tip},1}(t) \leq (1+\varepsilon) \, \sqrt{R_{\text{\rm tip},1}(t)}\] 
if $-t$ is sufficiently large (depending on $\varepsilon$). By Hamilton's Harnack inequality \cite{Hamilton2}, the function $t \mapsto R_{\text{\rm tip},1}(t)$ is monotone increasing. 
Consequently, 
\[d_{\text{\rm tip},1}((1+\varepsilon)t) - d_{\text{\rm tip},1}(t) \leq \varepsilon (1+\varepsilon) \, (-t) \, \sqrt{R_{\text{\rm tip},1}(t)}\] 
and 
\[d_{\text{\rm tip},1}(t) - d_{\text{\rm tip},1}((1-\varepsilon)t) \geq \varepsilon (1-\varepsilon) \, (-t) \, \sqrt{R_{\text{\rm tip},1}(t)}\] 
if $-t$ is sufficiently large. Using Proposition \ref{distance.to.tip}, we obtain 
\begin{align*} 
2\sqrt{1+\varepsilon} - 2 
&= \liminf_{t \to -\infty} \frac{d_{\text{\rm tip},1}((1+\varepsilon)t) - d_{\text{\rm tip},1}(t)}{\sqrt{(-t) \, \log(-t)}} \\ 
&\leq \varepsilon (1+\varepsilon) \liminf_{t \to -\infty} \sqrt{\frac{(-t) \, R_{\text{\rm tip},1}(t)}{\log (-t)}} 
\end{align*} 
and 
\begin{align*} 
2-2\sqrt{1-\varepsilon} 
&= \limsup_{t \to -\infty} \frac{d_{\text{\rm tip},1}(t) - d_{\text{\rm tip},1}((1-\varepsilon)t)}{\sqrt{(-t) \, \log(-t)}} \\ 
&\geq \varepsilon (1-\varepsilon) \limsup_{t \to -\infty} \sqrt{\frac{(-t) \, R_{\text{\rm tip},1}(t)}{\log (-t)}}.
\end{align*} 
Sending $\varepsilon \to 0$, we conclude that $\lim_{t \to -\infty} \sqrt{\frac{(-t) \, R_{\text{\rm tip},1}(t)}{\log (-t)}} = 1$. An analogous argument gives $\lim_{t \to -\infty} \sqrt{\frac{(-t) \, R_{\text{\rm tip},2}(t)}{\log (-t)}} = 1$. This completes the proof. 

\appendix 

\section{An elementary estimate for the one-dimensional heat equation}

\label{estimate.for.1d.heat.equation}

In this section, we prove an elementary estimate for the one-dimensional heat equation on the interval $[-1,1]$ with Dirichlet boundary conditions. This estimate is needed in the proof of Proposition \ref{key.estimate}. Recall that the heat kernel for the Dirichlet problem is given by 
\[K_t(x,y) = \frac{1}{\sqrt{4\pi t}} \bigg [ \sum_{k \in \mathbb{Z}} e^{-\frac{(x-y+4k)^2}{4t}} - \sum_{k \in \mathbb{Z}} e^{-\frac{(x+y+4k-2)^2}{4t}} \bigg ]\] 
for $x,y \in [-1,1]$ and $t > 0$. Note that $K_t(x,1) = K_t(x,-1) = 0$. We first record some basic properties of the Dirichlet heat kernel $K_t(x,y)$:

\begin{lemma}
\label{properties.of.Dirichlet.heat.kernel}
We can find a large constant $C$ such that the following holds:
\begin{itemize}
\item[(i)] $K_1(0,y) \geq \frac{1}{C} \, \cos \frac{\pi y}{2}$ for all $y \in (-1,1)$.
\item[(ii)] $\big | \frac{\partial^2}{\partial x^2} K_1(x,y) \big |_{x=0} \big | \leq C \, \cos \frac{\pi y}{2}$ for all $y \in (-1,1)$.
\item[(iii)] $-\frac{\partial}{\partial y} K_t(0,y) \big |_{y=1} \geq \frac{1}{C} \, t^{-\frac{3}{2}} \, e^{-\frac{1}{4t}}$ and $\frac{\partial}{\partial y} K_t(0,y) \big |_{y=-1} \geq \frac{1}{C} \, t^{-\frac{3}{2}} \, e^{-\frac{1}{4t}}$ for all $t \in (0,1]$.
\item[(iv)] $\big | \frac{\partial^2}{\partial x^2} \frac{\partial}{\partial y} K_t(x,y) \big |_{x=0, y=1} \big | \leq C \, t^{-\frac{7}{2}} \, e^{-\frac{1}{4t}}$ and $\big | \frac{\partial^2}{\partial x^2} \frac{\partial}{\partial y} K_t(x,y) \big |_{x=0, y=-1} \big | \leq C \, t^{-\frac{7}{2}} \, e^{-\frac{1}{4t}}$ for all $t \in (0,1]$.
\end{itemize}
\end{lemma}

\textbf{Proof.} 
We can find a small constant $\tau \in (0,1)$ such that the following holds:
\begin{itemize}
\item $-\frac{\partial}{\partial y} K_t(0,y) \big |_{y=1} \geq \frac{1}{C} \, t^{-\frac{3}{2}} \, e^{-\frac{1}{4t}}$ and $\frac{\partial}{\partial y} K_t(0,y) \big |_{y=-1} \geq \frac{1}{C} \, t^{-\frac{3}{2}} \, e^{-\frac{1}{4t}}$ for all $t \in (0,\tau]$.
\item $\big | \frac{\partial^2}{\partial x^2} \frac{\partial}{\partial y} K_t(x,y) \big |_{x=0, y=1} \big | \leq C \, t^{-\frac{7}{2}} \, e^{-\frac{1}{4t}}$ and $\big | \frac{\partial^2}{\partial x^2} \frac{\partial}{\partial y} K_t(x,y) \big |_{x=0, y=-1} \big | \leq C \, t^{-\frac{7}{2}} \, e^{-\frac{1}{4t}}$ for all $t \in (0,\tau]$.
\end{itemize}
In particular, $-\frac{\partial}{\partial y} K_\tau(0,y) \big |_{y=1}$ and $\frac{\partial}{\partial y} K_\tau(0,y) \big |_{y=-1}$ are positive numbers. Since $K_\tau(0,y) > 0$ for all $y \in (-1,1)$, we can find a small number $\varepsilon > 0$ such that $K_\tau(0,y) \geq \varepsilon \, \cos \frac{\pi y}{2}$ for all $y \in (-1,1)$. The maximum principle now implies $K_t(0,y) \geq \varepsilon \, e^{-\frac{\pi^2 t}{4}} \, \cos \frac{\pi y}{2}$ for all $y \in (-1,1)$ and all $t \in [\tau,1]$. In particular, $K_1(0,y) \geq \varepsilon \, e^{-\frac{\pi^2}{4}} \, \cos \frac{\pi y}{2}$. From this, statements (i), (iii), and (iv) follow.

To prove statement (ii), we observe that the function $y \mapsto \frac{\partial^2}{\partial x^2} K_1(x,y) \big |_{x=0}$ is smooth, and vanishes at $y=1$ and $y=-1$. Consequently, $\big | \frac{\partial^2}{\partial x^2} K_1(x,y) \big |_{x=0} \big | \leq C \, \cos \frac{\pi y}{2}$ for all $y \in (-1,1)$. This proves (ii). \\

\begin{proposition}
\label{1d.heat.equation}
Let $h(x,t)$ be a nonnegative solution of the heat equation $h_t - h_{xx} = 0$ on the rectangle $[-1,1] \times [-1,0]$. Then, for each $\mu \in (0,1)$, 
\[|h_{xx}(0,0)| \leq C\mu^{-2} \, h(0,0) + C e^{-\frac{1}{8\mu}} \, \sup_{\{-1,1\} \times [-1,0]} h,\] 
where $C$ is a constant.
\end{proposition}

\textbf{Proof.} 
The idea is to use the Greens representation formula for the Dirichlet problem for the heat equation on a rectangle. Fix a point $x \in (-1,1)$. For each $t \in (0,1]$, we define 
\[I(t) := \int_{-1}^1 K_t(x,y) \, h(y,-t) \, dy.\] 
Then 
\begin{align*} 
I'(t) 
&= \int_{-1}^1 \Big [ h(y,-t) \, \frac{\partial}{\partial t} K_t(x,y) + K_t(x,y) \, \frac{\partial}{\partial t} h(y,-t) \Big ] \, dy \\ 
&= \int_{-1}^1 \Big [ h(y,-t) \, \frac{\partial^2}{\partial y^2} K_t(x,y) - K_t(x,y) \, \frac{\partial^2}{\partial y^2} h(y,-t) \Big ] \, dy \\ 
&= h(1,-t) \, \frac{\partial}{\partial y} K_t(x,y) \Big |_{y=1} - h(-1,-t) \, \frac{\partial}{\partial y} K_t(x,y) \Big |_{y=-1}. 
\end{align*}
We now integrate this identity over $t \in (0,1]$. Since $\lim_{t \to 0} I(t) = h(x,0)$, it follows that 
\begin{align*} 
h(x,0) 
&= \int_{-1}^1 K_1(x,y) \, h(y,-1) \, dy \\ 
&- \int_0^1 h(1,-t) \, \frac{\partial}{\partial y} K_t(x,y) \Big |_{y=1} \, dt \\ 
&+ \int_0^1 h(-1,-t) \, \frac{\partial}{\partial y} K_t(x,y) \Big |_{y=-1} \, dt 
\end{align*} 
for $x \in (-1,1)$. We now put $x=0$. Using part (i) and (iii) of Lemma \ref{properties.of.Dirichlet.heat.kernel}, we obtain 
\begin{align*} 
h(0,0) 
&\geq \frac{1}{C} \int_{-1}^1 \cos \frac{\pi y}{2} \, h(y,-1) \, dy \\ 
&+ \frac{1}{C} \int_0^1 t^{-\frac{3}{2}} \, e^{-\frac{1}{4t}} \, [h(1,-t)+h(-1,-t)] \, dt. 
\end{align*} 
Similarly, using part (ii) and (iv) of Lemma \ref{properties.of.Dirichlet.heat.kernel}, we obtain 
\begin{align*} 
|h_{xx}(0,0)| 
&\leq C \int_{-1}^1 \cos \frac{\pi y}{2} \, h(y,-1) \, dy \\ 
&+ C \int_0^1 t^{-\frac{7}{2}} \, e^{-\frac{1}{4t}} \, [h(1,-t)+h(-1,-t)] \, dt. 
\end{align*} 
Putting these facts together, we conclude that 
\begin{align*} 
|h_{xx}(0,0)| 
&\leq C \int_{-1}^1 \cos \frac{\pi y}{2} \, h(y,-1) \, dy \\ 
&+ C \int_\mu^1 t^{-\frac{7}{2}} \, e^{-\frac{1}{4t}} \, [h(1,-t)+h(-1,-t)] \, dt \\ 
&+ C \int_0^\mu t^{-\frac{7}{2}} \, e^{-\frac{1}{4t}} \, [h(1,-t)+h(-1,-t)] \, dt \\ 
&\leq C \int_{-1}^1 \cos \frac{\pi y}{2} \, h(y,-1) \, dy \\ 
&+ C \mu^{-2} \int_\mu^1 t^{-\frac{3}{2}} \, e^{-\frac{1}{4t}} \, [h(1,-t)+h(-1,-t)] \, dt \\ 
&+ C \, e^{-\frac{1}{8\mu}} \int_0^\mu [h(1,-t)+h(-1,-t)] \, dt \\ 
&\leq C \mu^{-2} \, h(0,0) + C e^{-\frac{1}{8\mu}} \, \sup_{\{-1,1\} \times [-1,0]} h
\end{align*} 
for each $\mu \in (0,1)$. This completes the proof.

\section{The Bryant soliton}

\label{properties.of.Bryant.soliton}

In \cite{Bryant} Bryant showed that up to constant multiples, there is only one complete, steady, rotationally symmetric soliton in dimension three that is not flat. It has positive sectional curvature. The maximum scalar curvature is equal to $1$, and is attained at the center of rotation. The complete soliton can be written in the form $g = dz \otimes dz + B(z)^2 \, g_{S^2}$, where $z$ is the distance from the center of rotation. For large $z$, the metric has the following asymptotics: the aperature $B(z)$ has leading order term $\sqrt{2z}$, the orbital sectional curvature $K_{\text{\rm orb}}$ has leading order term $\frac{1}{2z}$, and the radial sectional curvature $K_{\text{\rm rad}}$ has leading order term $\frac{1}{4z^2}$.  

Sometimes it is more convenient to write the metric in the form $\Phi(r)^{-1} \, dr^2+ r^2 \, g_{S^2}$, where the function $\Phi(r)$ is defined by $\Phi(B(z)) = \big ( \frac{d}{dz} B(z) \big )^2$. The function $\Phi(r)$ is known to satisfy the equation
\[\Phi(r) \Phi''(r) - \frac{1}{2} \, \Phi'(r)^2 + r^{-2} \, (1-\Phi(r)) \, (r\Phi'(r)+2\Phi(r)) = 0.\] 
The orbital and radial sectional curvatures are given by $K_{\text{\rm orb}} = \frac{1}{r^2} \, (1-\Phi(r))$ and $K_{\text{\rm rad}} = -\frac{1}{2r} \, \Phi'(r)$. It is known that $\Phi(r)$ has the following asymptotics.  Near $r= 0$, $\Phi$ is smooth and has the asymptotic expansion
\[\Phi(r) = 1 + b_0 \, r^2 + o(r^2),\]
where $b_0$ is a negative constant. As, $r \to \infty$, $\Phi$ is smooth and has the asymptotic expansion
\[\Phi(r) = c_0 \, r^{-2} + 2 c_0^2 \, r^{-4} + o(r^{-4}),\]
where $c_0$ is a positive constant.

We will next find (for the convenience of the reader) the exact values of the constants $b_0$ and $c_0$ in the above asymptotics for the {\em Bryant soliton of maximal scalar curvature one}.  

Recall that the scalar curvature is given by $R = 2 K_{\text{\rm orb}} + 4 K_{\text{\rm rad}}$. The maximal scalar curvature is attained at $z = 0$, at which point $K_{\text{\rm orb}} = K_{\text{\rm rad}}$. The maximal scalar curvature being equal to $1$ is equivalent to $K_{\text{\rm orb}} = K_{\text{\rm rad}} = \frac{1}{6}$ at $z=0$. On the other hand, the asymptotic expansion of $\Phi(r)$ gives $K_{\text{\rm orb}} = \frac{1}{r^2} \, (1-\Phi(r)) = -b_0 + o(1)$ as $r \to 0$. Consequently, $b_0 = -\frac{1}{6}$.

Bryant's asymptotics imply that for $z$ sufficiently large, the aperture satisfies $r = (1+o(1)) \, \sqrt{2z}$, implying that $2z = (1+o(1)) \, r^2$. The radial sectional curvature satisfies $K_{\text{\rm rad}} = (1+o(1)) \, \frac{1}{4z^2} = (1+o(1)) \, r^{-4}$ for $r$ large.  On the other hand, the asymptotic expansion of $\Phi(r)$ implies $K_{\text{\rm rad}} = -\frac{1}{2r} \, \Phi'(r) = (1+o(1)) \, c_0 \, r^{-4}$ for $r$ large. Comparing the two formulae, we conclude that $c_0=1$.

Summarizing the above discussion we conclude the following asymptotics for the {\em Bryant soliton with maximal scalar curvature equal to one}:
\[\Phi(r) = \begin{cases} 1 - \frac{r^2}{6} + o(r^2) &\qquad \text{\rm as $r \to 0$,} \\ r^{-2} + 2r^{-4} + o(r^{-4}) &\qquad \text{\rm as $r \to \infty$.} \end{cases}\]

\end{document}